\documentclass[11pt]{article}
\usepackage{enumerate}
\usepackage{epsfig,alltt}
\usepackage[english]{babel}
\usepackage{blindtext}
\usepackage{dsfont}
\usepackage{amssymb}
\usepackage{amsfonts}
\usepackage{graphicx}
\usepackage{amsmath}
\usepackage{amsthm}
\usepackage{multirow}
\usepackage{natbib}
\usepackage{hyperref}
\usepackage[utf8]{inputenc}
\usepackage{amsthm,amsmath}
\usepackage[left=2.5cm,right=2.5cm,top=2.5cm,bottom=2.5cm]{geometry}

\usepackage{floatflt}
\usepackage{makeidx}
\usepackage{algorithm}

\usepackage{algorithmic}

\usepackage{url}
\usepackage{float} 
\theoremstyle{remark}

\usepackage{soul}

\usepackage{color}

\newcommand{\begeq}[1]{\begin{equation} \label{#1}}
	\newcommand{\fineq}{\end{equation}}

\title{Random effects estimation in a fractional diffusion model based on continuous observations}

\author{Nesrine CHEBLI$^{1,2}$\footnote{corresponding author}
	\footnote{e-mail address: nesrine.chebli@univ-poitiers.fr}, Hamdi FATHALLAH$^{2}$\footnote{e-mail address: hamdi.fathallah@essths.u-sousse.tn} and Yousri SLAOUI$^{1}$\footnote{e-mail address: Yousri.Slaoui@math.univ-poitiers.fr} \\
	$^{1}$ Laboratory of Mathematics and Applications, University of Poitiers, France\\
	$^{2}$ Laboratory of Mathematics Deterministic and Random Modeling, University of Sousse,\\ Tunisia
}

\numberwithin{equation}{section}
\begin{document}
	\newtheorem{theor}{\bf Theorem}
	\newtheorem*{theo}{\bf Theorem}
	\newtheorem{prop}{ \bf Proposition}
	\newtheorem{lem}{\bf Lemma}
	\newtheorem{coro}{\bf Corollary}
	\newtheorem{ex}{\bf Example}
	\newtheorem{prof}{\bf Proof}
	\newtheorem{defi}{\bf Definition}
	\newtheorem{rem}{\bf Remark}
	\newtheorem*{nota}{\bf Notations}
	\newtheorem*{assump}{\bf Assumptions}
	\newtheorem*{com}{\bf Comments on the assumptions}
	\newtheorem*{app}{\bf Appendix}
	\date{ }
	\maketitle
	
	\textit{\bf Abstract}: 
	The purpose of the present work is to construct estimators for the random effects in a fractional diffusion model using a hybrid estimation method where we combine parametric and nonparametric techniques. We precisely consider $n$ stochastic processes $\left\{X_t^j,\ 0\leq t\leq T\right\}$, $j=1,\ldots, n$ continuously observed over the time interval $[0,T]$, where the dynamics of each process are described by fractional stochastic differential equations with drifts depending on random effects. 
	We first construct a parametric estimator for random effects using maximum likelihood estimation techniques and study its asymptotic properties when the time horizon $T$ is sufficiently large. Then, on the basis of the obtained estimator for the random effects, we build a nonparametric estimator for their common unknown density function using Bernstein polynomials approximation. 
	Some asymptotic properties of the density estimator, such as its asymptotic bias, variance, and mean integrated squared error, are studied for an infinite time horizon $T$ and a fixed sample size $n$. The asymptotic normality of the estimator is established for a fixed $T$, a high frequency, and as long as the order of Bernstein polynomials is sufficiently large. We also investigate a non-asymptotic bound for the expected uniform error between the density function and its estimator. A numerical study is then presented in order to evaluate both qualitative and quantitative performance of the Bernstein estimator compared with the standard kernel estimator within and at boundaries of the support of the density function.\\
	
	\section{Introduction}
	\addcontentsline{toc}{section}{Introduction}
	Stochastic Differential Equations (SDEs) serve as a powerful mathematical tool to model dynamic systems that exhibit deterministic trends and random fluctuations. They can be thought of as a natural extension of ordinary differential equations that incorporate stochasticity into the modeling process, allowing us to account for random fluctuations and external influences. SDEs models are widely used in various applied problems, including physics, mathematical finance, biology, engineering, and economics, and the noise affecting the dynamics of the modeled systems is represented by a Brownian motion. 
	In the literature, statistical inference for stochastic processes modeled by SDEs, which we term as diffusion processes, has been thoroughly studied by now and summarized in several books; see, for example, \cite{bis08}, \cite{Kut04}, \cite{Pra99} and references cited therein. For these models, many authors have studied the drift parameter estimation in different models; see, for example, \cite{e}, \cite{f}, \cite{c}, \cite{d}, \cite{b}, \cite{a}.

	The significant constraint in employing stochastic diffusion models involving Brownian motion arises from the fact that it is characterized by the independence of its increments, resulting in uncorrelated random noise generated by this motion. 
	However, many phenomena emerging in a number of scientific domains may not have such a property. For example, in finance, the observed data present a long-term dependence in the sense that small variations are followed by small variations, and large variations are followed by large variations. 
	To integrate these properties in stochastic models, Mandelbrot and Van Ness proposed in \cite{man68} to modify the standard Brownian motion and popularized the normalized fractional Brownian motion (fBm) as a generalization of standard Brownian motion that exhibits long-range dependence. 
	The normalized fBm with Hurst index $H\in (0,1)$ is a centered Gaussian process with covariance function: $$\mathbb{E}\left(W_t^HW_s^H\right)=\dfrac{1}{2}\left(t^{2H}+s^{2H}-\lvert t-s \rvert^{2H} \right),\ \text{for all}\ t,s\geq 0.$$
	Therefore, SDEs driven by fBm are becoming more prominent in the statistical field as the most adequate models to describe this dependence, which justifies the significant interest in studying statistical inference problems for diffusion processes satisfying SDEs governed by fBm. \textcolor{blue}{Parametric estimation for these models using  the techniques of maximum likelihood estimation was investigated in \cite{ch13}, \cite{klep201}, \cite{Mis18} and references therein. In the nonparametric estimation for fractional SDEs, the reader can refer to \cite{C19}, \cite{marie23}, \cite{mishra11} and \cite{saus14} where the authors study Nadaraya-Watson-type estimators of the drift of a fractional SDE. See also \cite{Pr11} and the references therein for other estimation methods.}
	
	In another context, many studies are designed to examine changes over time in characteristics that are measured repeatedly for each study participant. For example, in a medical setting, data on neural potential, blood pressure measurements or cholesterol levels are obtained for each individual at different times and possibly under changing experimental conditions, making it difficult to determine whether the data are accurate, so that the probability distribution of the distribution of measurements has the same shape for each individual but the parameters of that distribution vary between individuals. To deal with such phenomena, SDEs with random effects have been introduced. 
	
	Statistical inference for random effects SDEs models has only recently garnered the attention of academics, compared to the huge literature on statistical inference for standard models. Because the variable of interest is not the observed variable, estimation problems for such models frequently occur. This might be caused, for example, by measurement errors during an experiment.
	The majority of the contributions with either discrete data or continuous data, assume a known model and focus on using parametric methods to estimate the parameters of the density of random effects, see for example the works of Delattre et al. \cite{Dal13}, Donnet and Samson \cite{Don08}, Genon-Catalot and Larédo \cite{GL14}, Picchini et al. \cite{Pic10} and Picchini and Ditlevsen \cite{PicD11}.\\
	The contributions to nonparametric estimation are fewer than those to parametric estimation and they address less broad models. 
	In this context, we call attention to \cite{C13} in which the authors developed nonparametric estimators for the density of random effects under restricted assumptions on the drift and diffusion coefficients and dealt with both additive and multiplicative random effects. In \cite{D14} and \cite{DG15}, Dion investigated a kernel estimator and a deconvolution estimator for the density of random effects in a stochastic diffusion model. 
	More recently, El Omari et al. studied in \cite{Om19} the properties of kernel and histogram estimators for the random effects density in a diffusion model governed by fBm and then based on their work, a more general random effects diffusion model described by SDEs driven by mixed fBm, was considered and studied in \cite{Pr(A2)0}.\\
	In the literature, kernel estimation is a common nonparametric estimation method. 
	However, kernel estimators present a main limitation which is the support problem when estimating functions with bounded support on at least one side. In fact, near the boundaries, the kernel is truncated, which leads to boundary bias and a significant underestimation of the true density values.  In order to overcome this problem, there has been a considerable development of methods for the estimation of density function. 
	One of the developed approaches is the use of Bernstein polynomials which were introduced in \cite{Ber12} in order to provide a probabilistic proof of the classical Weirstrass theorem which states that "Any continuous function on a segment $[a, b]$ is a uniform limit of polynomial functions on this segment". 
	Several publications have investigated nonparametric estimation based on Bernstein polynomials, such as 
	Babu et al. \cite{B02} and Leblanc \cite{L12} who studied the asymptotic properties of Bernstein estimators for density and distribution functions, \cite{BC06}, \cite{IK14} and \cite{J19} where the authors
	used Bernstein polynomials to construct a recursive density estimator and more recently the work of Slaoui \cite{S22} who proves moderate deviations principles for the recursive estimators of a distribution function defined by the stochastic approximation algorithm based on Bernstein polynomials. In all these contributions, it was shown that the Bernstein estimator has an interesting performance in the boundaries of the support of the distribution function or its density, in particular the absence of bias at the boundary points.
	
	To the best of our knowledge, nonparametric estimation based on Bernstein polynomials has not been investigated in the case of fractional SDEs with random effects yet.
	In the present work, we consider a general linear model described by the following fractional SDE with random effects \\
	\begin{equation}
		\label{**}
		\begin{aligned}
			dX_t &=S\left(\phi,t,X_t\right)dt + \sigma(t) dW^{H}_t,\ \ X_0=x_0,
		\end{aligned}
	\end{equation}
	where the process $\left\{X_t,\ 0\leq t\leq T\right\}$ is continuously observed on a time interval $[0,T]$,\\ $W^H=\left\{W_t^H,\ 0\leq t\leq T\right\}$ is a fBm with Hurst index $H\in (\frac{1}{2},1)$, $S: \Phi\subset\mathbb{R}^d\times [0,T]\times \mathbb{R}\longrightarrow \mathbb{R}$, called the drift coefficient, and $\sigma : [0,T]\longrightarrow \mathbb{R}_{+}$, called the diffusion coefficient, are known functions except $\phi$ which is an unobserved random effect with unknown density function $f$.\\
	Under local uniform Lipschitz continuity, linear growth and H\"{o}lder continuity conditions on the coefficients $S$ and $\sigma$, there exists a unique solution to the above fractional SDE, called a fractional diffusion process or simply a fractional diffusion (see \cite[Theorem 3.1.4 p.~201]{Mish08}).\\ 
	Our focus here is to construct an estimator for the density $f$ based on Bernstein polynomials and observations $\left\{X_t,\ 0\leq t\leq T\right\}.$
	
	This paper is organized as follows. In Section \ref{sec1} we introduce a particular case of the model ($\ref{**}$) where the function $S$ has a linear form and we list our notations and assumptions. In Section \ref{proc} we exhibit our main results, which we divide into three subsections. We begin by estimating the random effects and then constructing an estimate of their common density. The rest of the section is devoted to studying the asymptotic properties of the obtained estimator. 
	In Section \ref{s5} we illustrate the performance of our density estimator on simulated data and present a numerical comparison with the kernel estimator of the random effects density. Section \ref{conc} is devoted to the extensions and concluding remarks. To avoid interrupting the flow of this paper, all mathematical proofs are relegated to Section \ref{s6}.
	We close the paper with an Appendix, where we recall some classical limit theorems.
	
	\section{ Model, notations and assumptions}\label{sec1}

	In what follows, all random variables and processes are defined on a filtering probability space $(\Omega,\mathcal{F}, (\mathcal{F}_t ), \mathbb{P})$ satisfying the usual conditions and processes are $(\mathcal{F}_t )$-adapted. 
	On this space, we define $n$ stochastic processes $X^j=\left\{X^j_t,\ 0\leq t\leq T\right\}$, $j=1,\ldots,n$
	with dynamics ruled by the following fractional SDEs :
	\begin{equation}
		\begin{aligned}
			\begin{cases} 
				dX^{j}_t &=\left(a(X^j_t)+\phi_{j} b(t)\right)dt + \sigma(t) dW^{H,j}_t,\\
				X^{j}_{0} &= x^j \in \mathbb{R},
			\end{cases}
		\end{aligned}
		\label{eq1}
	\end{equation}
	where $\left(W^{H,j}\right)_{1\leq j \leq n}$, are $n$ independent normalized fBms with a common Hurst index $H\in (\frac{1}{2},1)$ and $\left(\phi_j\right)_{1\leq j \leq n}$ are $n$ unobserved independent and identically distributed (i.i.d.) real random variables, with a common density function $f$. Denote $F$ their distribution function. We assume that the processes $X^j$ are continuously observed on a time interval $[0, T ]$ with given $T > 0$. The sequences $\left(\phi_j\right)_{1\leq j \leq n}$ and $\left(W^{H,j}\right)_{1\leq j \leq n}$ are independent. The functions $b(.)$ and $a(.)$ are known on their own spaces and $\sigma(.)$ is a positive non-vanishing function. 
	When the observation time $t$ is fixed and due to the independence of the sequences $\left(\phi_j\right)_{1\leq j\leq n}$ and $\left(W^{H,j}\right)_{1\leq j\leq n}$, the random variables $X_t^j$ are i.i.d. We assume that $x^j=x_0$ so that the trajectories $\left\{X_t^j,\ 0\leq t\leq T\right\}$, $j=1, \ldots, n$ are i.i.d. 
	The main problem in estimating the common density $f$ arises from the fact that we do not observe the random effects $\phi_j$. The construction of an estimator of $\phi_j$ is therefore the initial step, after which we proceed to estimate their density function $f$.\\
	In the sequel, we introduce some needed notations and assumptions.
	\begin{nota}\ \\
		Throughout this paper, we mean by $o(.)$ and $O(.)$ the usual small-o and big-O, which means convergence in probability and stochastic boundedness, respectively. The notation $o_x(.)$ is used to mean that the limit depends on the point $x$. For any bounded function $h: [0,1]\rightarrow \mathbb{R}$, the norm \(\left\|h\right\|\) is defined by $\left\|h\right\|:= \underset{x\in [0,1]}{\sup} \left|h(x)\right|$. Finally, we denote by $\mathcal{N}$ and $\overset{\mathcal{L}}{\longrightarrow}$ the Gaussian distribution and the convergence in law, respectively. \\
		In order to construct the MLE of the random effects, we define for all $t\in [0,T]$,
		\begin{equation}\begin{aligned}
				\label{1}
				k_H(t,s)&=\kappa_H^{-1} s^{\frac{1}{2}-H}(t-s)^{\frac{1}{2}-H}\mathds{1}_{\left(0,t\right)}(s), \ \ \ \kappa_H&=2H \Gamma\left(\dfrac{3}{2} -H \right) \Gamma \left(H+ \frac{1}{2}\right),\end{aligned}
		\end{equation}
		\begin{equation}\begin{aligned}
				\label{2}
				w_t^{H}&=\lambda_H ^{-1} t^{2-2H}, \ \ \ \lambda_H&= \dfrac{2H \Gamma\left(3-2H\right)\Gamma\left(H+\frac{1}{2}\right)}{\Gamma\left(\dfrac{3}{2}-H\right)},
			\end{aligned}
		\end{equation}
		\begin{equation}
			\begin{aligned} 
				\label{4}
				J_1^j(t)&=\dfrac{d}{dw_t^H}\int_0^t k_H(t,s)\dfrac{a(X^j_s)}{\sigma(s)}ds,\ \ \ \ \ &
					J_2(t)&=\dfrac{d}{dw_t^H}\int_0^t k_H(t,s)\dfrac{b(s)}{\sigma(s)}ds.
				\end{aligned}
			\end{equation} Let $M^{H,j}=\left\{M^{H,j}_t,\ 0\leq t\leq T\right\}$ be a centered Gaussian process defined by 
			\begin{equation}
				\begin{aligned} 
					\label{3}
					M_0^{H,j}=0, \ M_t^{H,j}=\int_0^t k_H(t,s)dW_s^{H,j}.
				\end{aligned}
			\end{equation}
			\(M^{H,j}\) is a square-integrable martingale, called in \cite{nor19} the fundamental martingale or Molchan martingale, with quadratic variation given by $\langle M^{H,j} \rangle_t= w_t^H$. Furthermore, the natural filtration of martingale $M^{H,j}$ coincides with the natural filtration of fBm $W^{H,j}$.
		\end{nota}
		\begin{assump}\ \\
			\begin{enumerate}
				\item[$(A1)$] $\displaystyle \int_0^T J_2^2(t)dw_t^H <\infty$, a.s for any $T>0$.
				\item[$(A2)$] $\displaystyle \underset{T\to \infty}{\lim} \int^{T}_0 J_2^2(t)dw_t^H=\infty$ $a.s$.
				\item[$(A3)$] $f$ is twice continuously derivable on $[0,1]$.
			\end{enumerate}
		\end{assump}
		\begin{com}\ 
			\begin{itemize}
				\item Assumptions $(A1)$ and $(A2)$ are needed to estimate the random effects. Indeed, $(A1)$  ensures that the process $\displaystyle\int_0^t J_2(s)dM_s^{H,j}$, $t>0$, is a square integrable martingale and assumption $(A2)$ is needed to prove the consistency of the estimators of random effects.
				\item Assumption $(A3)$ is standard in the framework of nonparametric estimation of probability density using Bernstein polynomials.
			\end{itemize}
		\end{com}
		\section{ Hybrid estimation of the random effects density}\label{proc}
		This section is devoted to the construction of estimators for the random effects and their density function and to the study of the asymptotic behaviors of the obtained estimators.

		\subsection{Density approximation based on Bernstein polynomials}

		Assume that $Y_1,\ldots, Y_n$ are i.i.d random variables with distribution function $G$ and associated unknown density function $g$ supported on $[0,1]$.
		The ordinary Bernstein polynomial estimator of order $m$ for the density function $g$ is defined for all $x\in [0,1]$ by 
		\begin{equation*}
			\begin{aligned}
				\tilde{g}_{m,n}(x)=m\sum \limits_{k=0}^{m-1} \left[G_n\left(\dfrac{k+1}{m}\right)- G_n\left(\dfrac{k}{m}\right)\right] p_k\left(m-1,x\right),
			\end{aligned}
		\end{equation*}
		where $p_k\left(m,x\right)= C_m^k x^k \left(1-x\right)^{m-k}$ is the Bernstein polynomial and $G_n(y)=\dfrac{1}{n} \sum \limits_{j=1}^{n} \mathds{1}_{\left\{Y_j\leq y\right\}}$ is the empirical distribution function of $Y_j$. 
		For details of the properties of Bernstein polynomials, we refer the reader to \cite{Lor}. The estimator $\tilde{g}_{m,n}$ was introduced by Vitale \cite{V75} and was later studied by Babu et al. \cite{B02} and Leblanc \cite{L12} who developed results on its asymptotic properties.\\
		\textcolor{blue}{In the case where the random variable $Y$ is supported on the compact interval $[a, b]$, $a < b$, it can be easily transformed into a random variable $Z$ supported on $[0, 1]$ using the transformation $Z=\frac{1}{b-a}(Y-a)$. To cover random variables with values in $\mathbb{R}_{+}$ and $\mathbb{R}$ respectively, we use the transformations $Z=\frac{Y}{1+Y}$ and $Z=\frac{1}{2}+ \frac{1}{\pi} \arctan(Y)$. Once the random variable $Y$ is transformed, we can apply Bernstein polynomials to approximate its density function on $[0,1]$.}
		\textcolor{blue}{In this paper, we consider the case where the density of random effects $f$ is supported on $[0, 1]$. The Bernstein polynomials estimator of order $m > 0$ for $f$ is defined as follows }
		\begin{equation}
			\label{estb}
			\begin{aligned}
				\tilde{f}_{m,n}(x)=m\sum \limits_{k=0}^{m-1} \left[ F_n\left(\frac{k+1}{m}\right)- F_n\left(\frac{k}{m}\right)\right] p_k\left(m-1,x\right),
			\end{aligned}
		\end{equation}
		where $F_n(y)=\frac{1}{n} \sum \limits_{j=1}^{n} \mathds{1}_{\left\{\phi_j\leq y\right\}}$ is the empirical distribution function of $\phi_j$.
		Since the random effects are not observed, the estimator $\tilde{f}_{m,n}$ is not computable. Therefore, we first have to estimate $\phi_j$, $j=1,\ldots,n$ and then to estimate their density.
		\subsection{ MLE for the random effects and its asymptotic properties}
		In this subsection we construct for each $j=1, \ldots, n$ an approximation of the random variable $\phi_j$ using the techniques of maximum likelihood estimation. 
		The idea consists in considering for each fixed $j=1, \ldots, n$, the random variable $\phi_j$ as a deterministic unknown parameter $\varphi$ that will be estimated based on one observed trajectory $X^{j}$ which is the solution of the following fractional SDE 
		\begin{equation}
			dX^{j}_t=\left(a\left(X^{j}_t\right)+\varphi b\left(t\right)\right) dt + \sigma(t)dW_t^{H,j}, \ t\in [0,T],\ X_0^{j}=0,
		\end{equation}
		where $a(.)$ and $b(.)$ are as defined in the SDE (\ref{eq1}), $W^{H,j}=\left\{ W_t^{H,j}, \ t\in [0,T]\right\}$ is an fBm with Hurst parameter $H\in(\frac{1}{2},1)$ and $\sigma(.)$ is a non-vanishing positive function on $[0,\infty)$.\\ 
		Although fBm $W^{H,j}$ is not a semimartingale, by using an appropriate integral transformation, we can transform it into a martingale. \textcolor{blue}{Such integral transformation was introduced in \cite{nor19} and investigated in \cite{klep20}, \cite{klep202}, \cite{Pr11} and references therein.} 
		The obtained martingale is called in \cite{nor19} fundamental martingale and its natural filtration coincides with the natural filtration of fBm $W^{H,j}$. 
		Following the same techniques, we construct the martingale associated with fBm $W^{H,j}$, from which we then derive the MLE estimator of $\varphi$.\\
		Let for all $t\geq 0$,  
		$C_j(\varphi,t)=\displaystyle a\left(X^{j}_t\right)+\varphi b\left(t\right),$ 
		and assume that  $t\mapsto \dfrac{C_j(\varphi,t)}{\sigma(t)}$ is Lebesgue integrable on $[0,T]$ for any $T>0$ so that the process\begin{align*}
			Q^j_{H,\varphi}(t)=\dfrac{d}{dw_t^H}\int_0^t k_H(t,s) \dfrac{C_j(\varphi,s)}{\sigma(s)}ds, \ t\in [0,T],
		\end{align*}
		where $k_H(t, s)$ and $w_t^H$ are as defined in (\ref{1}) and (\ref{2}) respectively, is well defined. \\ Suppose the
		sample paths of the process $\left\{Q^j_{H,\varphi}(t), \ t\in [0,T]\right\}$ belong almost surely to $L^2\left([0, T], dw_t^H\right)$, i.e. $\displaystyle\int^T_0 \left(Q^j_{H,\varphi}(s)\right)^2dw_t^H < \infty$ a.s for all $T>0$ and define
		\begin{equation*}
			Z^j_t=\int_0^t \dfrac{k_H(t, s)}{\sigma(s)}dX^{j}_s, \ t\in [0,T].
		\end{equation*}
		The process $Z^j = \left\{Z^j_t, \ t\in [0,T]\right\}$ is a $(\mathcal{F}_t)$-semimartingale with the decomposition\begin{equation}
			\label{z}
			Z^j_t=\displaystyle\int_0^t Q^j_{H,\varphi}(s) dw_s^H + M_t^{H,j}= \int_0^t \left(J^{j}_1(s)+\varphi J_2(s)\right)dw_s^H + M_t^{H,j},
		\end{equation}
		where $J^{j}_1$ and $J^{j}_2$ are as defined in (\ref{4}) and $M^{H,j}$ is as defined in (\ref{3}).
		\begin{prop}\ 
			\label{pr1}
			\begin{enumerate}
				\item The MLE of the random effect $\phi_j$, $j=1, \ldots, n$, has the following form \begin{equation*}\label{phhi}
					\hat{\phi}_{j,T}=\dfrac{\displaystyle\int_0^T J_2(t) dZ^j_t-\int_0^T J_1^j(t)J_2(t) dw^H_t}{\displaystyle\int_0^T J^2_2(t) dw^H_t},\end{equation*}
				where $J_1^j(t)$ and $J_2(t)$ are as defined in (\ref{4}).
				\item  Under assumption $(A2)$, the estimator $\hat{\phi}_{j,T}$ is strongly consistent, that is \begin{equation*}\begin{aligned}\hat{\phi}_{j,T}\overset{\text{a.s}}{\longrightarrow} \phi_j\ \ \text{as}\ T\to \infty.\label{as}\end{aligned}\end{equation*}
				\item Suppose that there exists a deterministic function $h_t$ that satisfies 
				$$\underset{T\to \infty}{\lim} h_T  =0\ \text{a.s}\  \ \text{and}\ \ \underset{T\to \infty}{\lim} h_T^2 \int_0^T J_2^2(s) d w_s^H=c^2 <\infty,$$
				where $c$ is a positive constant.\\
				Then, under assumptions $(A1)$ and $(A2)$
				$$\displaystyle h_T^{-1}\left(\hat{\phi}_{j,T}-\phi_j\right)\ \overset{\mathcal{L}}{\longrightarrow} c\ \mathcal{G}\ \  \text{as}\  T\to \infty,$$
				where $\mathcal{G}$ is a standard normally distributed random variable.
			\end{enumerate}
		\end{prop}
		\begin{rem}\ \\
			In order to prove the third assertion of Proposition \ref{pr1} which states the asymptotic normality of the estimators $\hat{\phi}_{j,T}$, we use the classical Central Limit Theorem (CLT) for local continuous martingales. 
			We shall draw attention to the fact that we can obtain the asymptotic normality of $\hat{\phi}_{j,T}$ using the general version of CLT for martingales that was investigated by Touati in \cite{touati} where he used the characteristic function technicals instead of the classic Lindeberg condition (see Theorem \ref{touati1} in the Appendix). We can also use the CLT version established by Van Zanten in \cite{vanz}, Theorem 4.1.
	\end{rem}
	
	We now illustrate the results of Proposition \ref{pr1} in the following example.
	\begin{ex} {Fractional Vasicek model with random effects.}\ \\
		\label{examp}
		\begin{equation*}
			\begin{aligned}
				\begin{cases} 
					dX^{j}_t &=\left(-\beta X^j_t+\phi_{j} \right)dt + dW^{H,j}_t,\\
					X^{j}_{0} &= 0.
				\end{cases}
			\end{aligned}
			\label{ex}
		\end{equation*}
		where $W^{H,j}$ is a fBm with Hurst index $H> \frac{1}{2}$, $\beta$ is a known positive constant and $\phi_j$ is a random effect.\\
		From the results presented previously, the MLE of the random effect $\phi_j$ has the following form
		$$\hat{\phi}_{j,T}=\displaystyle\phi_j+\dfrac{\displaystyle\int_0^T J_2(t)dM_t^{H,j}}{\displaystyle\int_0^T J^2_2(t) dw^H_t},$$
		where, in this case,
		$$J_2(t)=\dfrac{d}{dw_t^H}\displaystyle\int_0^t k_H(t,s)ds,$$
		and $k_H(t,s)$ and $w_t^H$ are as defined in (\ref{1}) and (\ref{2}) respectively.\\
		However, $$\displaystyle \int_0^t k_H(t,s)ds=\kappa_H^{-1}\int_0^t  s^{\frac{1}{2}-H} (t-s)^{\frac{1}{2}-H}ds= \kappa_H^{-1}t^{2-2H} B\left(\dfrac{3}{2}-H,\dfrac{3}{2}-H\right)=\lambda_H^{-1}t^{2-2H}= w_t^H
		.$$
		Then \begin{align*}\displaystyle J_2(t)
			&=\dfrac{dt}{dw_t^H}\dfrac{d}{dt} \left(\lambda_H^{-1} t^{2-2H}\right)\\
			&=\dfrac{dt}{dw_t^H} \left((2-2H)\lambda_H^{-1}t^{1-2H}\right)\\
			&=\left((2-2H)\lambda_H^{-1}t^{1-2H}\right)^{-1}\left((2-2H)\lambda_H^{-1}t^{1-2H}\right)=1.
		\end{align*}
		Hence we obtain for $T>0$ and $j=1,\ldots,n$\begin{equation}
			\label{estim}
			\hat{\phi}_{j,T} =\displaystyle\phi_j+\dfrac{M_T^{H,j}}{w^H_T},\end{equation}
		where $w_T^H$ and  $M_t^{H,j}$ are as defined in (\ref{2}) and (\ref{3}).\\
		It is clear from (\ref{estim}) that for all $T>0$ and $j=1,\ldots, n$, 
		the quadratic risk $\mathbb{E}\left(\hat{\phi}_{j,T}-\phi_j\right)^2=\dfrac{\lambda_H}{T^{2-2H}}$ goes to $0$ as $T\to \infty$. 
		Since the process $M^{H,j}$ is a martingale with quadratic variation $\langle M^{H,j}\rangle_T=w_T^H$ and $H<1$, it follows that the function $\langle M^H\rangle_T$ tends to infinity as $T\to \infty$. By the strong law of large numbers for martingales, it follows that $$\dfrac{M_T^{H,j}}{w_T^H}=\dfrac{M_T^{H,j}}{\langle M^H\rangle_T} \longrightarrow 0,\ \ \ \text{a.s}\ \ \text{as}\ \ T\to \infty.$$
		Hence, $\hat{\phi}_{j,T} \overset{\text{a.s}}{\longrightarrow} \phi_j$, as $T\to \infty$.\\
		Recall that $M^{H,j}$ is a centered Gaussian process with variance function $w^H_T$, it follows that $$\dfrac{M_T^{H,j}}{\sqrt{w_T^H}} \overset{d}{\longrightarrow}\mathcal{N}(0,1).$$
		As a result
		$$T^{1-H}\left(\hat{\phi}_{j,T}-\phi_j\right)=\dfrac{T^{H-1}}{\sqrt{w_T^H}}\dfrac{M_T^{H,j}}{ \sqrt{w_T^H}}=\sqrt{\lambda_H}\dfrac{M_T^{H,j}}{ \sqrt{w_T^H}}\overset{d}{\longrightarrow} \mathcal{N}(0,\lambda_H),$$
		which proves the asymptotic normality of $\hat{\phi}_{j,T}$.
	\end{ex}
	
	The next subsection is devoted to Bernstein density estimators of random effects $\phi_j$.
	
	\subsection{Density estimator and its asymptotic properties
	} \label{s3}
	The estimator of the common density $f$ of the random effects is obtained by injecting $\hat{\phi}_{j,T}$ into the expression $(\ref{estb})$. 
	Let for all $x\in [0,1]$,   
	\begin{equation} 
		\label{estb1}
		\begin{aligned}
			\hat{f}_{m,n}(x)
			= m\sum \limits_{k=0}^{m-1}\left[\hat{F}_n\left(\dfrac{k+1}{m}\right)- \hat{F}_n\left(\dfrac{k}{m}\right)\right]p_k\left(m-1,x\right),
		\end{aligned}
	\end{equation}
	where $\hat{F}_n(x)=\dfrac{1}{n} \sum \limits_{j=1}^{n}\mathds{1}_{\left\{\hat{\phi}_{j,T}\leq x\right\}}$ is the empirical estimator for the cumulative distribution of $\hat{\phi}_{j,T}$. 
	The remainder of this section is devoted to studying the asymptotic properties of the constructed estimator (\ref{estb1}). 
	We first begin by studying the bias and variance of $\hat{f}_{m,n}$ which are stated below.
	\begin{prop}
		\label{prop1} \ \\
		Let assumptions $(A1)-(A3)$ hold.\begin{description}
			\item[1)]	For $x\in [0,1]$, we obtain \begin{equation*}
				\label{Bia} 
				\begin{aligned}
					\lim\limits_{T\to \infty}	Bias\left(\hat{f}_{m,n}(x)\right)= m^{-1}\dfrac{\left(1-2x\right)}{2}f'(x)+ o\left(m^{-1}\right).
				\end{aligned}	
			\end{equation*}	
			\item[2)] \begin{enumerate}
				\item[i)] For $x\in (0,1)$
				, we obtain
				\begin{equation*} 
					\begin{aligned}
						\lim\limits_{T\to \infty}	Var\left(\hat{f}_{m,n}(x)\right)=  m^{\frac{1}{2}}n^{-1}\gamma(x) +o_x\left(m^{\frac{1}{2}}n^{-1}\right),
					\end{aligned}
				\end{equation*}
				where $\gamma(x)=f(x)\psi(x)$ with $\psi(x)=\left(4\pi x (1-x)\right)^{-\frac{1}{2}}.$
				\item[ii)] For $x\in \left\{0,1\right\}$
				, we obtain
				\begin{equation*}
					\begin{aligned}
						\lim\limits_{T\to \infty}	Var\left(\hat{f}_{m,n}(x)\right)= mn^{-1} f(x)+O(mn^{-1}).
					\end{aligned}
				\end{equation*}
			\end{enumerate}	
			
		\end{description}
	\end{prop}
	
	A simple way to characterize the global error of an estimator of an unknown density is to study its Mean Integrated Squared Error (MISE) as presented in the following corollary.
	\begin{coro}\ \\
		\label{cor1}
		$ $ 
		Under assumptions $(A1)-(A3)$, we obtain 
		\begin{equation*}
			\begin{aligned}
				\lim\limits_{T\to \infty}	MISE\left(\hat{f}_{m,n}\right)= m^{\frac{1}{2}}n^{-1} C_1 + m^{-2} C_2+ o\left(m^{\frac{1}{2}}n^{-1}\right) + o\left(m^{-2}\right), 
			\end{aligned}
		\end{equation*}
		where $C_1=\displaystyle \int_0^1\left(\dfrac{1-2x}{2}\right)^2 f(x)^2 dx$ and $C_2=\displaystyle\int_0^1 f(x)\psi(x) dx$.
	\end{coro}
	\begin{rem}\ \\
		If $C_1>0$, the optimal choice of $m$ based on the minimization of the MISE, is
		$$m_{\text{opt}}=\left(\dfrac{4C_2}{C_1}\right)^{\frac{2}{5}}n^{\frac{2}{5}}.$$
		The corresponding MISE is given as follows
		$$\lim\limits_{T\to \infty}MISE\left(\hat{f}_{m_{\text{opt}},n}\right)=\dfrac{5}{4} 4^{\frac{1}{5}} C_1^{\frac{4}{5}} C_2^{\frac{1}{5}} n^{-\frac{4}{5}}+ o\left(n^{-\frac{4}{5}}\right).$$
		Therefore the estimator $\hat{f}_{m,n}$ achieves the optimal rate of convergence in terms of MISE for density functions which is of order of $n^{-\frac{4}{5}}$ when $m$ is chosen proportional to $n^{\frac{2}{5}}$. Compared with the kernel density estimator, EL Omari et al. demonstrated in \cite{Om19} that the convergence rate for the kernel estimator of the random effects density is of the order of $n^{-\frac{2\beta}{2\beta +1}}$ when the bandwidth $h$ is chosen proportional to $n^{-\frac{1}{2\beta +1}}$ and the kernel is chosen of order equal to $\lfloor\beta \rfloor$, where $\lfloor\beta \rfloor$ denotes the greatest integer strictly less than the real number $\beta$. If we let $\beta=2$ and $h=m^{-1}$, we conclude that the optimal bandwidth for $\hat{f}_{m,n}$ is $h=n^{-\frac{2}{5}}$ instead of $h=n^{-\frac{1}{5}}$ for kernel estimators.
	\end{rem}
	In the remainder of this section, we assume that $m$ is a function of $n$ such that $m=m_n \to \infty$ as $n\to\infty$ and study the way in which $\hat{f}_{m,n}$ behaves with respect to $f$. To that end, we consider the following additional assumption: 
	\begin{enumerate}
		\item [$(A4)$] the random variable $\hat{\phi}_{1,T}$ has a positive density function $g$ that is continuous at $x$ and note $G$ its probability distribution function. 
	\end{enumerate}
	\textcolor{blue}{
		\begin{theor}\ \\
			\label{Teo1}
			Under assumptions $(A1)-(A4)$ and as $m,n\longrightarrow \infty$ in such a way that $mn^{-1}\longrightarrow 0$, the following assertions hold.
			\begin{enumerate}
				\item For $x\in(0,1)$ 
				\begin{eqnarray}
					n^{\frac{1}{2}}m^{-\frac{1}{4}}\left[\hat{f}_{m,n}(x)- \mathbb{E}\left(\hat{f}_{m,n}(x)\right)\right] \overset{\mathcal{L}}{\longrightarrow} \mathcal{N}\left(0,g(x)\psi(x)\right). 
				\end{eqnarray}
				\item For $x\in\left\{0,1\right\}$ 
				\begin{eqnarray}
					n^{\frac{1}{2}}m^{-\frac{1}{2}}\left[\hat{f}_{m,n}(x)- \mathbb{E}\left(\hat{f}_{m,n}(x)\right)\right] \overset{\mathcal{L}}{\longrightarrow} \mathcal{N}\left(0,g(x)\right). 
				\end{eqnarray}
			\end{enumerate}
	\end{theor}}
	\textcolor{blue}{The asymptotic normality of $\hat{f}_{m,n}$ is established in the following corollary, which follows as a deduction from Theorem \ref{Teo1}.}
	\textcolor{blue}{
		\begin{coro}\ \\
			\label{co2}
			Under assumptions $(A1)-(A4)$, the following assertions hold as $m,n\to\infty$ such that $mn^{-1}\to 0$.\begin{enumerate}
				\item If $n^{-\frac{1}{2}} m^{-\frac{1}{4}} \longrightarrow c$ for some constant $c \geq 0$, then for $x\in\left(0,1\right)$
				\begin{align*}
					n^{\frac{1}{2}} m^{-\frac{1}{4}} \left[\hat{f}_{m,n}(x)-f(x)\right] \overset{\mathcal{L}}{\longrightarrow} \mathcal{N}\left(\delta,g(x)\psi(x)\right).
				\end{align*}
				\item If $n^{-\frac{1}{2}} m^{-\frac{1}{2}} \longrightarrow c$ for some constant $c \geq 0$, then for $x\in\left\{0,1\right\}$\begin{align*}
					n^{\frac{1}{2}} m^{-\frac{1}{2}} \left[\hat{f}_{m,n}(x)-f(x)\right] \overset{\mathcal{L}}{\longrightarrow} \mathcal{N}\left(\delta,g(x)\right),
				\end{align*}
			\end{enumerate}	
			where $\delta=c\left(g(x)-f(x)\right)$. 
	\end{coro}}
	\textcolor{blue}{In the following theorem, we present a non-asymptotic risk bound for the expected uniform error between $\hat{f}_{m,n}$ and the true density function $f$.}
	\textcolor{blue}{\begin{theor}\ \\
			\label{t3}
			Let assumptions $(A1)-(A3)$ hold and assume that for all $s\in[0,t]$ there exists a positive constant $C$ such that $ C\leq \dfrac{b(s)}{\sigma(s)}$. Then 
			\[     \mathbb{E} \left\|\hat{f}_{m,n}-f \right\|\leq  \frac{\lambda_H m^4}{C^2 T^{2 - 2H}}+ m^{\frac{3}{2}}n^{-1/2} +m^{-1}\left(\dfrac{\left\|f'\right\|}{2}+\dfrac{\left\|f''\right\|}{8}\right)+ O(m^{-1}).  \]
	\end{theor}}
	\section{Numerical simulation}
	\label{s5}
	In the following section, we present a numerical study of the Bernstein estimator of the common density of random effects in the fractional Vasicek model previously discussed in Example \ref{examp}. 
	Our first task is to compute the MLE $\hat{\phi}_{j,T}$ after simulating the true random effects. 
	To that end, we assume that the trajectories $\left\{X^j_t,\ 0\leq t \leq T \right\}$ are observed at time points $t_k=k\dfrac{T}{n}$, $k=1,\ldots,n+1$ simultaneously, with different sample sizes $n$ and fixed $T=100$. The fBm is simulated using the Yuima package, and the Molchan martingale is approximated using a simple Riemann sum.\\
	For the common density $f$, we investigate four distributions with different behavior: \begin{itemize}
		\item Beta density $\mathcal{B}(1,2)$.
		\item Beta density $\mathcal{B}(3,5)$.
		\item Beta mixture density $0.5\mathcal{B}(3,9)+0.5\mathcal{B}(9,3)$.
		\item Truncated normal mixture density 0.6*$\mathcal{N}(0.5, 0.1)$ + 0.4* $\mathcal{N}(0.9,0.03)$.
	\end{itemize}
	\textcolor{blue}{For each density function, the Bernstein estimator $\hat{f}_{m,n}$ is implemented with different sample sizes $n$ and an optimal data-driven choice of Bernstein polynomials order $m$. To determine this optimal $m$, we used the least squares cross-validation method (LSCV). For our density estimator \(\hat{f}_{m,n}(x)\), the LSCV criterion is given by
		\[
		\text{LSCV}(m) = \int_0^1 \hat{f}_{m,n}(x)^2 dx - \dfrac{2}{n} \sum_{j=1}^n \hat{f}_{m,n}^{(-j)}(\hat{\phi}_{j,T}),
		\]
		where\begin{itemize}
			\item \(\hat{f}_{m,n}(x)\) is the Bernstein estimator that uses all the data.
			\item\(\hat{f}_{m,n}^{(-j)}(x)\) is the Bernstein estimator computed without the \(j\)-th observation.
			\item The term \(\dfrac{2}{n} \sum\limits_{j=1}^n \hat{f}_{m,n}^{(-j)}(\hat{\phi}_{j,T})\) is a cross-validation term that penalizes overfitting.
		\end{itemize}
		The optimal \(m\) is the one that minimizes \(\text{LSCV}(m)\). Let \(m=\underset{m}{arg\min}\ \text{LSCV}(m).\) \\
		For each example of the four densities, Bernstein estimator is compared with the standard kernel density estimator defined by \begin{align*}
			\hat{f}_h(x)=\dfrac{1}{nh}\sum \limits_{i=1}^{n} K\left(\dfrac{x-\hat{\phi}_{i,T}}{h}\right)
		\end{align*}
		where $K$ is the Gaussian kernel and the bandwidth $h$ is selected using Silverman's rule $$h_{\text{opt}}=1.06\ \min\left(sd(\hat{\phi}_{j,T}),\ 1.34\ IQR(\hat{\phi}_{j,T})\right)\ n^{-\frac{1}{5}},$$ 
		where $IQR$ is the Interquartile Range and $sd$ is the standard deviation.}\\
	In what follows, we begin with presenting a qualitative comparison between the two estimators for each example of density function.
	
	\begin{figure}[H]
		\centering
		\includegraphics[scale=0.56]{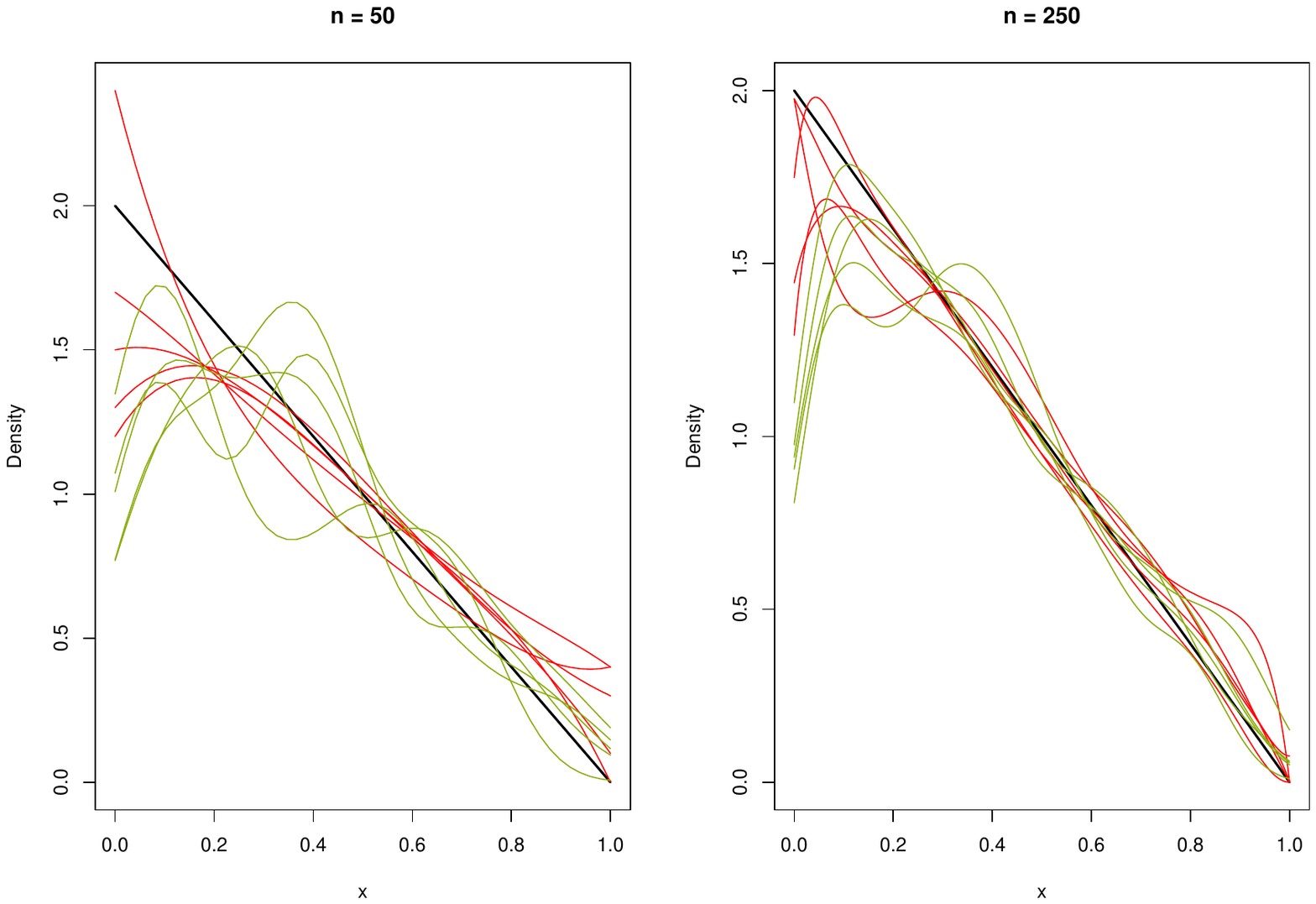}
		\caption{ Qualitative comparison between the Bernstein estimator $\hat{f}_{m,n}$ and the Kernel estimator $\hat{f}_h$ of the density Beta(1,2). The black line represents the true density $f(x)$, the red line represents the Bernstein estimator $\hat{f}_{m,n}(x)$ and the green line represents the Kernel Density Estimator $\hat{f}_h(x)$}
		\label{fig1}
	\end{figure}
	
	\begin{figure}[H]
		\centering
		\includegraphics[scale=0.56]{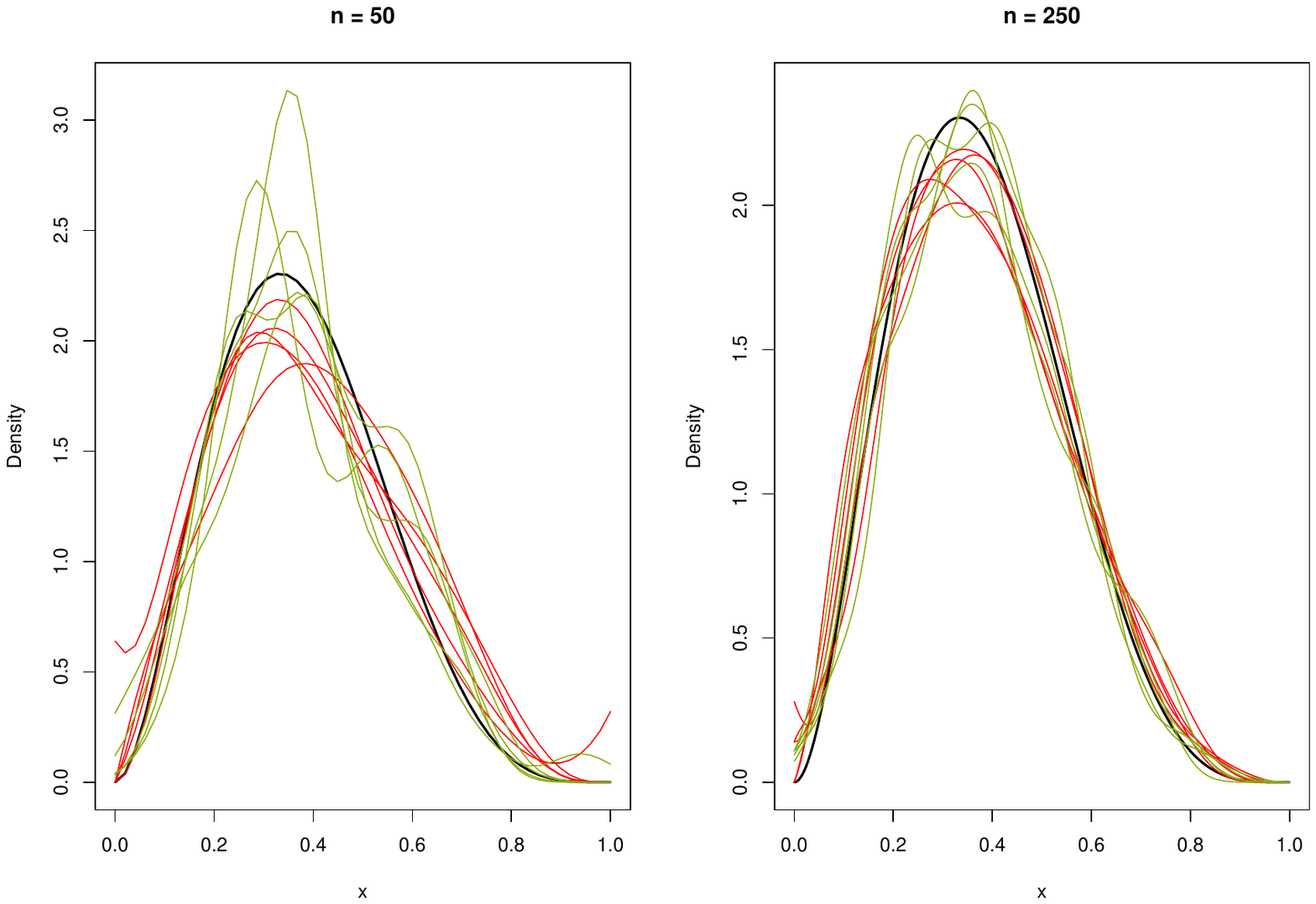}
		\caption{Qualitative comparison between the Bernstein estimator $\hat{f}_{m,n}$ and the Kernel estimator $\hat{f}_h$ of the density Beta(3,5). The black line represents the true density $f(x)$, the red line represents the Bernstein estimator $\hat{f}_{m,n}(x)$ and the green line represents the Kernel Density Estimator $\hat{f}_h(x)$}
		\label{fig2}
	\end{figure}
	
	\begin{figure}[H]
		\centering
		\includegraphics[scale=0.56]{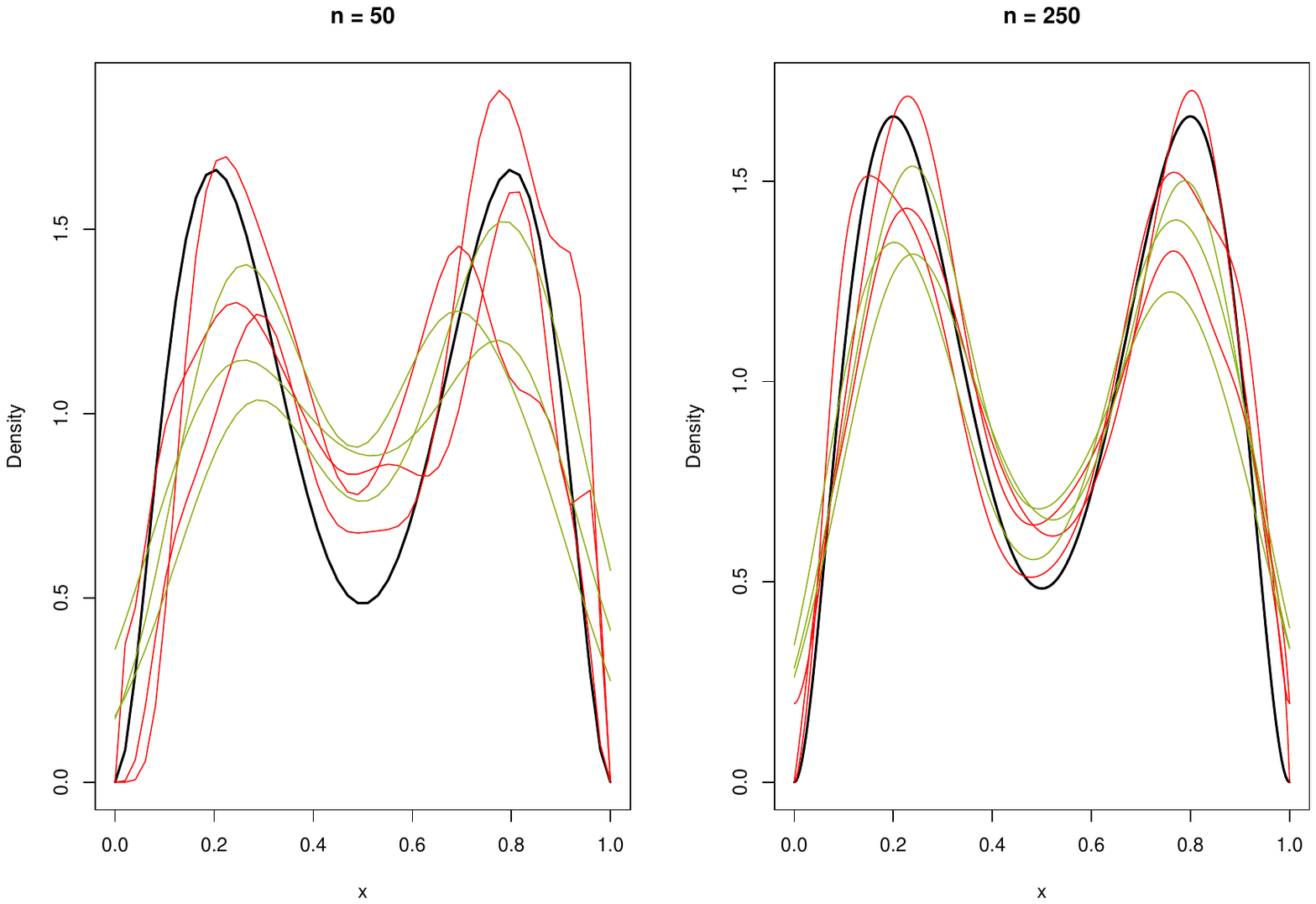}
		\caption{Qualitative comparison between the Bernstein estimator $\hat{f}_{m,n}$ and the Kernel estimator $\hat{f}_h$ of the density 0.5*Beta(3,9) + 0.5* Beta(9,3). The black line represents the true density $f(x)$, the red line represents the Bernstein estimator $\hat{f}_{m,n}(x)$ and the green line represents the Kernel Density Estimator $\hat{f}_h(x)$.}
		\label{fig3}
	\end{figure}
	\begin{figure}[H]
		\centering 
		\includegraphics[scale=0.56]{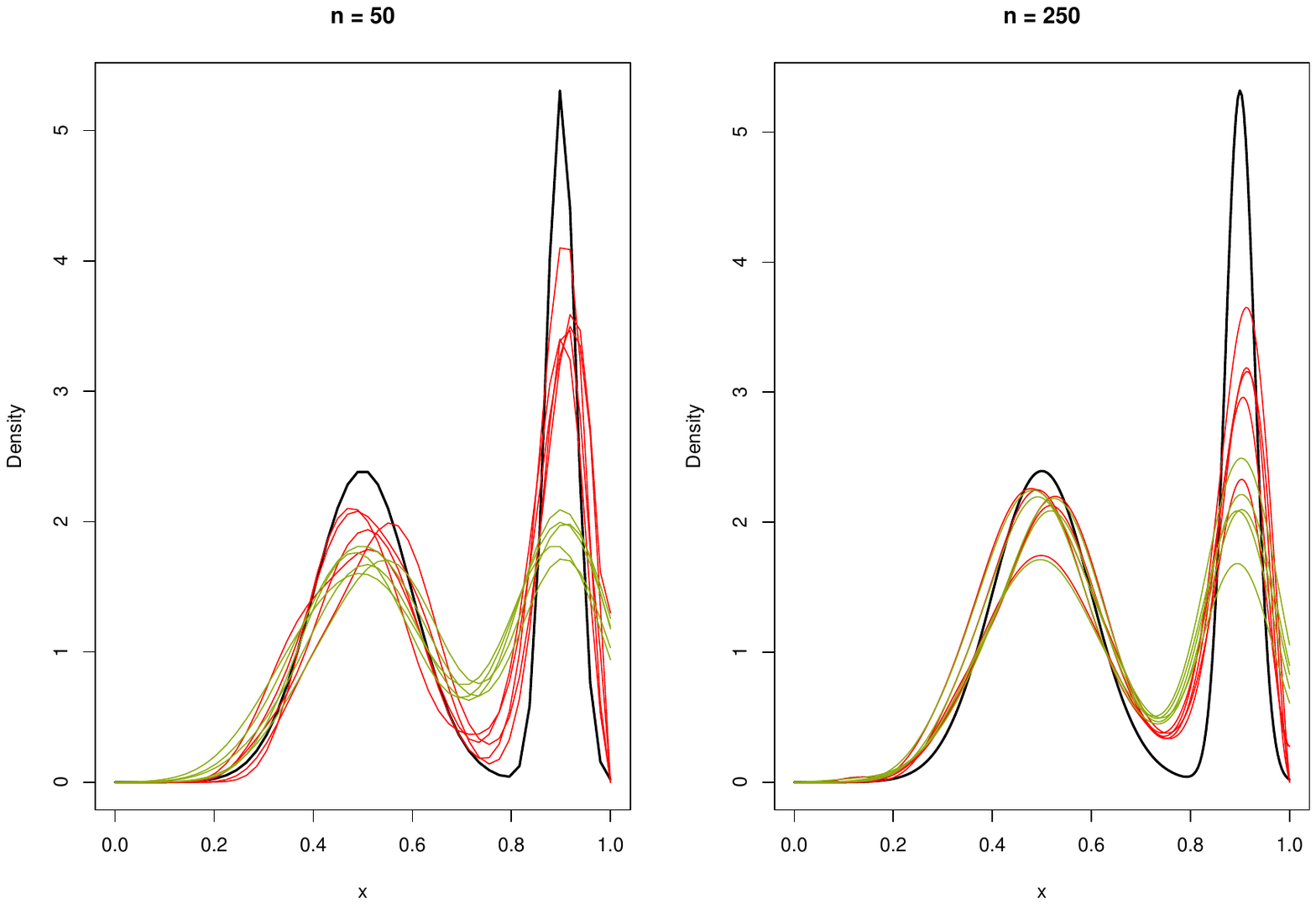}
		\caption{Qualitative comparison between the Bernstein estimator $\hat{f}_{m,n}$ and the Kernel estimator $\hat{f}_h$ of the density 0.6*$\mathcal{N}(0.5, (0.1)^2)$ + 0.4* $\mathcal{N}(0.9,(0.03)^2)$. The black line represents the true density $f(x)$, the red line represents the Bernstein estimator $\hat{f}_{m,n}(x)$ and the green line represents the Kernel Density Estimator $\hat{f}_h(x)$}
		\label{fig4}
	\end{figure}
	\paragraph{\textcolor{blue}{Comments on the Figures:}}\
	\textcolor{blue}{\begin{itemize}
			\item  From Figure \ref{fig1}, we observe that the Bernstein estimator $\hat{f}_{m,n}$ (red) consistently follows the general shape of the true density more smoothly than the kernel estimator $\hat{f}_{h}$ (green), which exhibits noticeable variability and oscillations. At $x=0$, it's obvious that the Bernstein estimator better approximates the peak value, highlighting its improved boundary performance, whereas the Kernel estimator shows more variability and less accuracy.
			\item In Figure \ref{fig2}, the Bernstein estimator is smoother and better follows the overall shape of the true density for both sample sizes $n=50$ and $n=250$, but slightly underestimates the peak that represents the Beta(3,5) density around $x\approx 0.4$, whereas the kernel estimator overestimates the peak of the density for $n=50$ and still exhibit more fluctuations and instability despite its improvement for greater sample size $n=250$.
			\item For the Beta mixture density function presented in Figure \ref{fig3}, the Bernstein estimator overestimates the true density and the kernel estimator underestimates it for a small sample size $n=50$. As the sample size increases, both estimators show significant improvements especially the Bernstein estimator, which aligns closely with the two peaks and valley of the density. At boundaries, it's obvious that the kernel estimator still struggle to approximate the boundary values of the density.
			\item From Figure \ref{fig4}, we conclude that the kernel estimator captures sharp features better but at the cost of high variance, whereas the Bernstein estimator provides a smoother and more stable estimate while slightly losing precision in the peaks. Near the boundaries, $\hat{f}_{m,n}$ provides more stable estimates than $\hat{f}_{h}$ and attains the true values of the density at the boundary points.
			\item As the sample size $n$ increases from 50 to 250, the optimal order of Bernstein polynomials $m$ increases and we obtain closer estimates to the true densities.  
			\item In all cases, the Bernstein estimator shows faster convergence near the boundaries as $n$ increases.
			\item In all cases, the Bernstein estimator performs much better near the boundaries, while the kernel estimator tends to underestimate the density functions and exhibits boundary bias. A detailed numerical evaluation of the performance of the two estimators near the boundaries is presented in Table \ref{bound} below. 
	\end{itemize}}
	\textcolor{blue}{From a quantitative perspective, the performance of the two density estimators is evaluated using three different error metrics: 
		Integrated Squared Error (ISE), to evaluate the global accuracy of the estimators over [0,1], Mean Squared Error (MSE) and Mean Absolute Error (MAE) to measure the average error across all evaluation points. 
		Each error metric is calculated for different sample sizes $n\in\left\{50, 200, 500\right\}$ across the four different densities as shown in the following table.}
	\begin{table}[H]
		\centering
		\textcolor{blue}{  \begin{tabular}{|c|c|c|c|c|c|c|c|}
				\hline \multicolumn{2}{|c|}{}& \multicolumn{2}{c|}{ISE}&\multicolumn{2}{c|}{MSE} &\multicolumn{2}{c|}{MAE}\\
				\hline density &$n$& $\hat{f}_{m,n}$& $\hat{f}_h$&$\hat{f}_{m,n}$& $\hat{f}_h$&$\hat{f}_{m,n}$& $\hat{f}_h$\\
				\hline$(1) $&50&\textbf{0.068906}&0.086823&\textbf{0.029395}&0.079103&\textbf{0.113015}& 0.182032\\
				&200&\textbf{0.032575}&0.058011&\textbf{0.046186}&0.074050&0.164423&0.162584\\
				&500&\textbf{0.019774}&0.040545&\textbf{0.015067}&0.030554&\textbf{0.059808}&0.093648\\ 
				\hline$(2)$&50&\textbf{0.3030793}& 0.3244993&0.027516&\textbf{0.010081}&0.143378&\textbf{0.080494}\\
				&200&\textbf{0.296797}& 0.329376&\textbf{0.004424}& 0.008330&\textbf{0.057569}&0.069698\\
				&500&\textbf{0.274272}&0.311323&\textbf{0.024601}&0.034552&\textbf{0.068032}&0.097142\\  
				\hline$(3)$&50&\textbf{0.095608}&0.109241&\textbf{0.060238}&0.077975&\textbf{0.189213}&0.238338\\
				&200&\textbf{0.049940 }&0.057823&\textbf{0.029470}&0.046315&\textbf{0.134674}&0.187355\\
				&500&\textbf{0.034432}& 0.031541&\textbf{0.448455}&0.570038&\textbf{0.312530}&0.367121\\
				\hline$(4)$&50& \textbf{0.5782835}& 0.7817266&\textbf{0.515275}&0.705154&\textbf{0.364487}&0.483081\\
				&200&\textbf{0.5349293}&0.5631679&\textbf{ 0.567284}&0.543785&\textbf{0.398433}&0.436049\\
				&500&0.5106915&\textbf{0.4413549}&0.457656&\textbf{0.395983}&0.361263&\textbf{0.315452}\\
				\hline
		\end{tabular}}
		\textcolor{blue}{\caption{The average integrated squared error (ISE), the mean squared error (MSE) and the mean absolute error (MAE) of the Bernstein estimator $\hat{f}_{m_{opt},n}$ and the kernel estimator $\hat{f}_h$ for each sample size $n$ and for optimal $m$ and $h$.}
			\label{Tab:1}}
	\end{table}
	\textcolor{blue}{From Table \ref{Tab:1}, we can conclude that \begin{itemize}
			\item Among the three error metrics (ISE, MSE, and MAE), the Bernstein estimator $\hat{f}_{m,n}$ outperforms the kernel estimator and consistently shows lower values as indicated by bold values in the table.
			\item Both estimators show improvement with larger sample size $n$ but the kernel estimator still underperforms the Bernstein estimator.
			\item For both estimators, ISE, MSE, and MAE decrease as the sample size increases. 
	\end{itemize}}
	
	\textcolor{blue}{In order to thoroughly evaluate the performance of both estimators $\hat{f}_{m,n}$ and $\hat{f}_{h}$ near the boundaries, we present in the following tables a numerical comparison between the values of density functions at the exact boundary points $x=0$ and $x=1$, and at points near the boundaries $x=0.01$ and $x=0.99$.}
	\begin{table}[H]
		\centering
		\resizebox{\textwidth}{!}{%
			\textcolor{blue}{\begin{tabular}{|c|c|cccc|cccc|cccc|}
					\hline
					\multicolumn{2}{|c|}{} & \multicolumn{4}{c|}{\textbf{True Density}} & \multicolumn{4}{c|}{\textbf{Bernstein Estimator}} & \multicolumn{4}{c|}{\textbf{Kernel Estimator}} \\
					\hline
					density & $n$ & $x = 0$ & $x = 0.01$ & $x = 0.99$ & $x = 1$ & $x = 0$ & $x = 0.01$ & $x = 0.99$ & $x = 1$ & $x = 0$ & $x = 0.01$ & $x = 0.99$ & $x = 1$ \\
					\hline  
					\multirow{3}{*}{(1)} & 50  & 2 & 1.98 & 0.02  & 0  & \bf{2.4}  & \bf{2.333246} & \bf{0.396531} & 0.4  & 1.3468  & 1.420620 & 0.131555 & \bf{0.1167}  \\
					& 100 & 2 & 1.98 & 0.02  & 0  & \bf{1.98}  & \bf{1.75784} & \bf{0.026977} & \bf{0}  & 0.9641  & 0.953755 & 0.042668 & 0.0784  \\
					& 250 & 2 & 1.98 & 0.02  & 0  & \bf{1.792} & \bf{1.504638} & \bf{0.097609} & \bf{0} & 0.8494  & 0.889556 & 0.054277 & 0.0338  \\
					\hline  
					\multirow{3}{*}{(2)} & 50  & 0 & 0.010086 & 0.000001  & 0  & \bf{0}  & \bf{0.009542} & \bf{0.000016} & \bf{0}   & 0.041350  &  0.055893 &   0.000130 & 0.000071  \\
					& 100 & 0 & 0.010086 & 0.000001  & 0  & 0.31 & 0.310079 & \bf{0.0000001} & \bf{0}  & \bf{0.153994}  & 0.180859 & 0.000062 & \bf{0.000029} \\
					& 250 & 0 & 0.010086 & 0.000001  & 0  & 0.64 & 0.159405 & \bf{0.000029} & 0 & \bf{0.095671}  & \bf{0.128181}  & 0.000034 & 0.000015  \\
					\hline 
					\multirow{3}{*}{(3)} & 50  & 0 & 0.022837 & 0.022837  & 0  & \bf{0}  & \bf{0.463748} & \bf{0.00121} & \bf{0} & 0.5177  & 0.553786 & 0.614688 & 0.5585  \\
					& 100 & 0 & 0.022837 & 0.022837  & 0  & \bf{0.13}  & \bf{0.293611} & \bf{0.385897} & \bf{0.13}  & 0.3497  & 0.375827 & 0.540479 & 0.3765  \\
					& 250 & 0 & 0.022837 & 0.022837  & 0  & \bf{0.13}  & \bf{0.209428} & \bf{0.266990} & \bf{0.13}  & 0.3064  & 0.353690 & 0.332978 & 0.2891  \\
					\hline    
					\multirow{3}{*}{(4)} & 50  & 0.00009 & 0.000015 & 0.059106 & 0.02057  & \bf{0}     & \bf{0} & \bf{0.176836} & \bf{0}     & 0.00096 &  0.001388 & 1.052810 & 0.9403 \\
					& 100 & 0.00009 & 0.000015 & 0.059106  & 0.02057  & \bf{0}     & \bf{0} & \bf{0.344660} & \bf{0}     & 0.00094 & 0.000054 & 1.128598 & 0.8538 \\
					& 250 & 0.00009 & 0.000015 & 0.059106  & 0.02057  & \bf{0}     & \bf{0}  & \bf{0.267674} & \bf{0}     & 0.00112 & 0.001669 & 0.869418 & 0.7141 \\
					\hline   
				\end{tabular}
		}}
		\textcolor{blue}{\caption{Comparison of True densities, Bernstein estimators, and Kernel estimators at boundaries}
			\label{bound}}
	\end{table}
	\textcolor{blue}{From Table \ref{bound}, we can deduce that the Bernstein estimator provides a better approximation and generally outperforms the kernel estimator near and at boundary points across all sample sizes $n$, while the kernel estimator tends to underestimate the true density values at these points. Therefore, the kernel estimator still struggle with accurately approximating the true density near the boundaries, which is one of its noticeable limitations that highlights the advantage of the Bernstein estimator in boundary correction.}
	\section{Conclusion and perspective} \label{conc}
	In this work we discuss the problem of random effects estimation in a linear fractional stochastic diffusion model. We particularly consider $n$ stochastic processes $\left\{X_t^j,\ t\in [0,T]\right\}$, $j=1,\ldots, n$ observed continuously on a time interval $[0,T]$ and described by linear stochastic differential equations governed by fBm and including random effects $\phi_j$, $j=1,\ldots,n$ in the drift coefficients.
	We first construct for each $j$ a parametric estimator for the random effect $\phi_j$ using maximum likelihood estimation method. Then using a plug-in technique we build a nonparametric estimator based on Bernstein polynomials for the common density of the random effects. The asymptotic behavior of both parametric and nonparametric estimators is also investigated. In order to illustrate our theoretical results and ensure how good is our Bernstein density estimator, we present a numerical study based on a qualitative and quantitative comparison with the kernel density estimator. According to the obtained numerical results, Bernstein estimator is the most efficient in the context of this paper. 
	It demonstrates higher accuracy in approximating the true density function specifically at the boundaries of its compact support.  \\
	An extension of the present work to diffusion models governed by a mixed fractional Brownian motion with random effects where a simultaneous parametric estimation of all model parameters is investigated and using Tchebychev polynomials to approximate the law of random effects in the basis of discrete observations is an ongoing work. Another future research direction would be to extend our study to the case of censored data (see, e.g., Slaoui \cite{sl22}).
	\section{Proofs}\label{s6}
	This section is devoted to the detailed proofs of our results.
	\subsection{Proof of Proposition \ref{pr1}}
	\begin{enumerate}
		\item  Let $\mathbb{P}_{\varphi}^T$ be the measure induced by the process $\left\{X^j_t,\ 0\leq t \leq T\right\}$ when $\varphi$ is the true parameter. Hence, the Radon-Nikodym derivative of $\mathbb{P}_{\varphi}^T$ with respect to $\mathbb{P}_{0}^T$ is given by\begin{equation*}
			L_T\left(\varphi\right):=\dfrac{d\mathbb{P}_{\varphi}^T}{d\mathbb{P}_{0}^T}=\exp\left[\int_0^T Q^j_{H,\varphi}(s) dZ^j_s-\dfrac{1}{2}\int_0^T Q^j_{H,\varphi}(s)^2dw^H_s\right].
		\end{equation*}
		Therefore \begin{align*}
			\log\left(L_T(\varphi)\right)=\int_0^T \left(J^{j}_1(t)+\varphi J_2(t)\right)dZ_t - \dfrac{1}{2}\int_0^T\left(J^{j}_1(t)+\varphi J_2(t)\right)^2 dw_t^H,
		\end{align*}
		and the likelihood equation is given as follows
		\begin{align*}
			\int_0^T J_2(t) dZ_t - \int_0^T  \left(J^j_1(t)+\varphi J_2(t)\right)J_2(t) dw_t^H =0.
		\end{align*}
		At last the MLE $\hat{\varphi}_T$ of $\varphi$ is given for $T>0$ by
		\begin{equation*}
			\begin{aligned}
				\label{fi}
				\hat{\varphi}_T=\dfrac{\int_0^T J_2(t) dZ^j_t-\int_0^T J^{j}_1(t)J^{\varphi}_2(t) dw^H_t}{\int_0^T J^2_2(t) dw^H_t}.
			\end{aligned}
		\end{equation*}
		
		\item From the decomposition ($\ref{z}$) and the expression of $\hat{\phi}_{j,T}$ given in the first assertion of Proposition \ref{pr1},
		the error term is given for $j=0, \ldots, n$ and $T>0$ by 
		\begin{equation}\begin{aligned} \label{er}\hat{\phi}_{j,T}-\phi_j=\dfrac{\int_0^T J_2(t)dM_t^{H,j}}{\int_0^T J^2_2(t) dw^H_t}.\end{aligned}\end{equation}
		Since the process $M^{H,j}=\left\{M^{H,j}_t,\ t\in [0,T]\right\}$ is a centered square integrable local martingale with quadratic variation $\langle M^{H,j}\rangle_t=w^H_t$,
		and $\left\{J_2(t),\ t\in [0,T]\right\}$ is a measurable process that satisfies assumptions $(A1)$ and $(A2)$, 
		it follows from the strong law of large numbers for the continuous martingale $M^{H,j}$ (see Theorem \ref{Liptser and Shiryaev} in the Appendix) that $$\dfrac{\displaystyle\int_0^T J_2(t)dM_t^{H,j}}{\displaystyle\int_0^T J^2_2(t) dw^H_t}
		\overset{\text{a.s}}{\longrightarrow} 0 \quad  \text{as} \quad T\to \infty,$$
		which implies the consistency of the estimator $\hat{\phi}_{j,T}$ and we obtain the second assertion of proposition \ref{pr1}.
		\item To prove the third assertion of Proposition \ref{pr1} that presents the asymptotic normality of the estimator $\hat{\phi}_{j,T}$, we consider the error term previously expressed in (\ref{er}) and we assume that we have found a function $h_t$ that verifies the conditions of the proposition, namely
		$$\underset{T\to \infty}{\lim} h_T  =0\ \text{a.s}\  \ \text{and}\ \ \underset{T\to \infty}{\lim} h_T^2 \int_0^T J_2^2(s) d w_s^H=c^2 <\infty.$$
		We then apply the Central Limit Theorem (CLT) for the local continuous martingale $M^H$, which completes the proof of Proposition \ref{pr1}.
	\end{enumerate}

	\subsection{Proof of Proposition \ref{prop1} }
	In order to simplify notations, we denote for $m\geq 1$, $k=0,\ldots, m-1$ and $j=1,\ldots, n$  \begin{equation*}\begin{aligned}
			\mathcal{B}_k=\left]\dfrac{k}{m}, \dfrac{k+1}{m}\right] \ \ \ \text{and}\ \ \ 
			Z^T_{j,m}(x)=\sum \limits_{k=0}^{m-1}\mathds{1}_{\left\{\hat{\phi}_{j,T}\in \mathcal{B}_k\right\}}p_k\left(m-1,x\right).
		\end{aligned}
	\end{equation*}
	Therefore the density estimator $\hat{f}_{m,n}$ defined in (\ref{estb1}) can be written, for all $x\in [0,1]$ and $m\geq 1$, as follows
	\begin{align*}
		\hat{f}_{m,n}(x)=\dfrac{m}{n}\sum \limits_{j=1}^{n}Z_{j,m}^T(x).
	\end{align*}
	\begin{enumerate}
		\item For the first assertion of our Proposition, its sufficient to study the asymptotic behavior of the mean of $\hat{f}_{m,n}$ as $T$ tends to $\infty$.\\
		Let for all $x\in [0,1]$, $\mathbb{E}\left(\hat{f}_{m,n}(x)\right)=\dfrac{m}{n}\sum \limits_{j=1}^{n} \mathbb{E}\left(Z^T_{j,m}(x)\right)=m\mathbb{E}\left(Z^T_{1,m}(x)\right).$\\
	However
	$$\mathbb{E}\left(Z^T_{1,m}(x)\right)=\sum \limits_{k=0}^{m-1} \mathbb{E}\left(\mathds{1}_{\left\{\hat{\phi}_{1,T}\in \mathcal{B}_k\right\}}\right)p_k\left(m-1,x\right)=\sum \limits_{k=0}^{m-1}\mathbb{P}\left(\hat{\phi}_{1,T}\in \mathcal{B}_k \right)p_k\left(m-1,x\right).$$
	According to the consistency of $\hat{\phi}_{j,T}$ given by the second assertion of Proposition \ref{pr1}, we have for all $j=1,\ldots, n$,\ \  $\hat{\phi}_{j,T} \overset{\mathcal{L}}{\longrightarrow} \phi_j, \ \text{as}\ \ T\longrightarrow\infty$. Since $\mathbb{P}\left(\phi_j \in \delta \mathcal{B}_k\right)=0$, the Portmanteau theorem implies that \begin{equation}
		\label{conv}
		\underset{T\to \infty}{\lim} \mathbb{P}\left(\hat{\phi}_{j,T}\in \mathcal{B}_k\right)=\mathbb{P}\left(\phi_{j}\in \mathcal{B}_k\right).\end{equation}
	Consequently, $\underset{T\to \infty}{\lim} \mathbb{E}\left(Z^T_{1,m}(x)\right)=\sum \limits_{k=0}^{m-1}\mathbb{P}\left(\phi_{1}\in \mathcal{B}_k\right)p_k\left(m-1,x\right)$.\\
	It follows that \begin{equation} \label{lex}\begin{split}
			\underset{T\to \infty}{\lim}\mathbb{E}\left(\hat{f}_{m,n}(x)\right)= m\underset{T\to \infty}{\lim} \mathbb{E}\left(Z^T_{1,m}(x)\right)&=m\sum \limits_{k=0}^{m-1}\mathbb{P}\left(\phi_{1}\in \mathcal{B}_k\right)p_k\left(m-1,x\right)\\&=m\sum \limits_{k=0}^{m-1}\left(F\left(\dfrac{k+1}{m}\right)-F\left(\dfrac{k}{m}\right)\right)p_k\left(m-1,x\right).\end{split}\end{equation}
	Now, our goal is to study the asymptotic behaviour of the last equality.
	Using Taylor-Young's theorem for the distribution function $F$, we get for all $0\leq k \leq m-1$
	\begin{equation}
		\label{d1}
		F\left(\dfrac{k+1}{m}\right)-F\left(\dfrac{k}{m}\right)=\dfrac{1}{m}f\left(\dfrac{k}{m}\right)+\dfrac{1}{2m^2}f'\left(\dfrac{k}{m}\right)+\textbf{o} \left(m^{-2}\right).\end{equation}
	Again thanks to Taylor-Young's formula applied to the density function $f$ and its derivative $f'$ we obtain respectively
	\begin{equation}
		\label{d2}
		f\left(\dfrac{k}{m}\right)=f(x)+\left(\dfrac{k}{m}-x\right)f'(x)+o\left(\dfrac{k}{m}-x\right).\end{equation}
	and \begin{equation}\label{d3}f'\left(\dfrac{k}{m}\right)=f'(x)+o\left(\dfrac{k}{m}-x\right).\end{equation}
	Combining (\ref{d1}), (\ref{d2}) and (\ref{d3}), we obtain 
	\begin{align*}
		m\left(F\left(\dfrac{k+1}{m}\right)-F\left(\dfrac{k}{m}\right)\right)&=f(x)+\left(\dfrac{k}{m}-x\right)f'(x)+ \dfrac{1}{2m}f'(x)+o\left(\dfrac{k}{m}-x\right)\\&+ \dfrac{1}{2m}o\left(\dfrac{k}{m}-x\right)+o(m^{-1}).
	\end{align*}
	Substituting these expansions into equation (\ref{lex}) allows us to write
	\begin{align*}\underset{T\to \infty}{\lim}\mathbb{E}\left(\hat{f}_{m,n}(x)\right)&=f(x)\sum \limits_{k=0}^{m-1}p_k\left(m-1,x\right)+f'(x)\sum \limits_{k=0}^{m-1}\left(\dfrac{k}{m}-x\right)p_k\left(m-1,x\right)\\ &+o\left(\sum \limits_{k=0}^{m-1}\left(\dfrac{k}{m}-x\right)p_k\left(m-1,x\right)\right)+\dfrac{1}{2m} f'(x)\sum \limits_{k=0}^{m-1}p_k\left(m-1,x\right)\\ &+o\left(\dfrac{1}{2m}\sum \limits_{k=0}^{m-1}p_k\left(m-1,x\right)\right).\end{align*}
	Now using the fact that  $ \sum \limits_{k=0}^{m-1}p_k\left(m-1,x\right)=1$ and $\sum \limits_{k=0}^{m-1}\left(\dfrac{k}{m}-x\right)p_k\left(m-1,x\right)=-\dfrac{x}{m}$, we obtain
	\begin{align*}
		\underset{T\to \infty}{\lim}\mathbb{E}\left(\hat{f}_{m,n}(x)\right)=f(x) + \left(\dfrac{1-2x}{2m}\right)f'(x)+o\left(-\dfrac{x}{m}\right)+o(m^{-1}),
	\end{align*}
	which leads to the needed expression of the bias.
	\item For the first item of the second assertion of Proposition \ref{prop1}, we have, for all $x\in (0,1)$,
	\begin{equation}
		\label{var}
		\begin{split}
			Var\left(\hat{f}_{m,n}(x)\right)=\dfrac{m^2}{n^2}Var\left(\sum \limits_{j=1}^{n}Z^T_{j,m}(x)\right)= \dfrac{m^2}{n}\left(\mathbb{E}\left(Z^T_{1,m}(x)^2\right)- \left(\mathbb{E}\left(Z^T_{1,m}(x)\right)\right)^2\right).  
		\end{split}
	\end{equation}
	Moreover\begin{align*} \mathbb{E}\left(Z^T_{1,m}(x)^2\right)&=\mathbb{E}\left[\left(\sum \limits_{k=0}^{m-1}\mathds{1}_{\left\{\hat{\phi}_{1,T}\in \mathcal{B}_k\right\}}p_k\left(m-1,x\right)\right)^2\right]\\
		&=\mathbb{E}\left[\sum \limits_{k=0}^{m-1}\mathds{1}_{\left\{\hat{\phi}_{1,T}\in \mathcal{B}_k\right\}}p_k^2\left(m-1,x\right)\right]\\
		&=\sum \limits_{k=0}^{m-1} \mathbb{P}\left(\hat{\phi}_{1,T}\in \mathcal{B}_k\right)p_k^2\left(m-1,x\right).
	\end{align*}
	Using again (\ref{conv}), it follows that \begin{align*}\lim\limits_{T\to \infty}\mathbb{E}\left(Z^T_{1,m}(x)^2\right)&=\sum \limits_{k=0}^{m-1} \mathbb{P}\left(\phi_1\in \mathcal{B}_k\right)p_k^2\left(m-1,x\right)\\
		&=\sum \limits_{k=0}^{m-1}\left[ F\left(\dfrac{k+1}{m}\right)-F\left( \dfrac{k}{m}\right)\right]p_k^2\left(m-1,x\right).\end{align*}
	An application of Taylor-Young’s formula implies that
	\begin{equation}
		\label{lim}
		\begin{aligned}
			\lim\limits_{T\to \infty}\mathbb{E}\left(\left(Z^T_{1,m}(x)\right)^2\right)&=\dfrac{1}{m} f(x) \sum \limits_{k=0}^{m-1}p_k^2\left(m-1,x\right) + O\left(\sum \limits_{k=0}^{m-1} \left|\dfrac{k}{m}-x\right|p_k^2\left(m-1,x\right)\right)\\ &+O\left(\dfrac{1}{m}\sum \limits_{k=0}^{m-1}p_k^2\left(m-1,x\right)\right).
		\end{aligned}
	\end{equation}
	For the first and the last terms of (\ref{lim}), we use Lemma 3.1 in \cite{B02} which states that for all $x\in(0,1)$ $$\sum \limits_{k=0}^{m-1}p_k^2\left(m-1,x\right)=m^{-\frac{1}{2}} \psi(x)\left(1+o(1)\right), $$
	where $\psi(x)=\left(4 \pi x(1-x)\right)^{-\frac{1}{2}}$.\\
	For the second term of (\ref{lim}),
	the Cauchy-Schwarz inequality combined with the fact that for all $k=0,\ldots,m-1$, $0\leq p_k\left(m-1,x\right)\leq 1$ implies that
	\begin{align*}
		\left(\sum \limits_{k=0}^{m-1} \left|\dfrac{k}{m}-x\right|p_k^2\left(m-1,x\right)\right)^2&\leq \sum \limits_{k=0}^{m-1} p_k^3\left(m-1,x\right) \sum \limits_{k=0}^{m-1} p_k\left(m-1,x\right) \left(\dfrac{k}{m}-x\right)^2\\
		&= O\left(m^{-1} \sum \limits_{k=0}^{m-1} p_k^3\left(m-1,x\right)\right)
		= O\left(m^{-\frac{3}{2}}\right),
	\end{align*}
	which leads to $$\sum \limits_{k=0}^{m-1} \left|\dfrac{k}{m}-x\right|p_k^2\left(m-1,x\right)=O\left(m^{-\frac{3}{4}} \right).$$
	Then \begin{equation}\label{z1}\begin{split} \lim\limits_{T\to \infty}\mathbb{E}\left(\left(Z^T_{1,m}(x)\right)^2\right)&
			=m^{-\frac{3}{2}}f(x)\psi(x)\left(1+o(1)\right)+O(m^{-\frac{7}{4}})+O_x(m^{-\frac{3}{2}}).\end{split}\end{equation}
	To complete the proof, we develop the second term on the left hand side of (\ref{var}) as follows
	\begin{equation}\label{z2}\begin{aligned}
			\lim\limits_{T\to \infty}\left[\mathbb{E}\left(Z^T_{1,m}(x)\right)\right]^2=\left[\dfrac{1}{m} f(x)+O(m^{-1})\right]^2 =O(m^{-2}).
		\end{aligned}
	\end{equation}
	Combining (\ref{z1}) with (\ref{z2}), we obtain , for all $ x\in(0,1)$
	\begin{align*}
		\lim\limits_{T\to \infty}Var\left(\hat{f}_{m,n}(x)\right)&= \lim\limits_{T\to \infty}m^2n^{-1} Var\left(Z^T_{1,m}(x)\right)\\ &=m^{\frac{1}{2}}n^{-1}f(x)\psi(x) +o_x\left(m^{\frac{1}{2}}n^{-1}\right).
	\end{align*}
	
	Now for the second item of the assertion which concerns all $x\in \left\{0,1\right\}$, we have 
	\begin{align*}
		\sum \limits_{k=0}^{m-1}p_k^2\left(m-1,0\right)=\sum \limits_{k=0}^{m-1}p_k^2\left(m-1,1\right)=1,
	\end{align*}
	and thanks to (\ref{lim}), we obtain\begin{align*}
		\lim\limits_{T\to \infty}\mathbb{E}\left(Z^T_{1,m}(x)^2\right)=m^{-1} f(x)+O(m^{-1}).
	\end{align*} 
	Consequently, for $x\in \left\{0,1\right\}$
	\begin{align*}
		\lim\limits_{T\to \infty}Var\left(\hat{f}_{m,n}(x)\right)=mn^{-1} f(x)+O(mn^{-1}),
	\end{align*}
	which completes the proof of Proposition \ref{prop1}.
\end{enumerate}
\subsection{Proof of corollary \ref{cor1}}
Using Fubini's Theorem and the results of Proposition \ref{prop1}, we obtain
\begin{equation}
	\label{mise1}
	\begin{aligned}
		\lim\limits_{T\to \infty}MISE\left[\hat{f}_{m,n}\right]=	\lim\limits_{T\to \infty}\int_0^1 \left[Var\left(\hat{f}_{m,n}(x)\right) +  Bias^2\left(\hat{f}_{m,n}(x)\right)\right] dx.
	\end{aligned}
\end{equation}
For the bias term, we have
\begin{equation}
	\begin{aligned}
		\label{bi}
		\int_0^1 Bias^2\left(\hat{f}_{m,n}(x)\right) dx&= m^{-2} \int_0^1\left(\dfrac{1-2x}{2}\right)^2 f^2(x) dx +o(m^{-2})\\
		&= m^{-2}C_1+ o\left(m^{-2}\right), 
	\end{aligned}
\end{equation}
where $C_1=\displaystyle \int_0^1\left(\dfrac{1-2x}{2}\right)^2 f^2(x) dx$,
and for the variance term, we have
\begin{equation}
	\label{va}
	\begin{aligned}
		\lim\limits_{T\to \infty}\int_0^1 Var\left(\hat{f}_{m,n}(x)\right) dx &=m^2n^{-1}\int_0^1 Var\left(Z^T_{1,m}(x)\right) dx \\ 
		&= m^{\frac{1}{2}}n^{-1}\int_0^1\left[\gamma(x) \left(1+o(1)\right) +O\left(m^{\frac{1}{2}}n^{-1}\right)\right]dx\\
		&= m^{\frac{1}{2}}n^{-1}\left[C_2+o(1)\right]
	\end{aligned}
\end{equation}
where $C_2=\displaystyle\int_0^1 f(x)\psi(x) dx$.\\
By combining (\ref{mise1}), (\ref{bi}) and (\ref{va}), we obtain the announced result.
\textcolor{blue}{\subsection{Proof of Theorem \ref{Teo1}}}\
\textcolor{blue}{For all $x\in [0,1]$, we have
	\begin{align*}
		\hat{f}_{n,m}(x)-\mathbb{E}\left(\hat{f}_{n,m}(x)\right) &
		=mn^{-1} \sum \limits_{j=1}^{n} Z^T_{j,m}(x)-m\mathbb{E}\left(Z^T_{1,m}(x)\right)\\
		&=mn^{-1} \sum \limits_{j=1}^{n} T_{j,m},
	\end{align*}
	where $T_{j,m}=Z^T_{j,m}(x)-\mathbb{E}\left(Z^T_{1,m}(x)\right).$}\ 
\textcolor{blue}{It follows that \begin{align*}
		n^{\frac{1}{2}}m^{-1} \left(\hat{f}_{n,m}(x)-\mathbb{E}\left(\hat{f}_{n,m}(x)\right)\right)=\sum \limits_{j=1}^{n} \dfrac{T_{j,m}}{n^{\frac{1}{2}}}.
	\end{align*}
	In order to apply Lindeberg Feller's Central Limit Theorem, we shall first verify the Lindeberg condition.\\
	To that end, we define $Y_{j,m}:=\dfrac{T_{j,m}}{n^{\frac{1}{2}}}$, which are i.i.d. random variables with mean 0 and denote $s_n^2=\sum \limits_{j=1}^{n} \mathbb{E}\left(Y_{j,m}^2\right)$.}\\
\textcolor{blue}{It comes that $$n^{\frac{1}{2}}m^{-1} \left(\hat{f}_{n,m}(x)-\mathbb{E}\left(\hat{f}_{n,m}(x)\right)\right)=\sum \limits_{j=1}^{n} Y_{j,m}.$$
	To verify the Lindeberg condition, it is enough to show that for all $\varepsilon >0$,
	\begin{equation}
		\label{lin}
		\begin{aligned}
			\dfrac{1}{s_n^2}\sum \limits_{j=1}^{n} \mathbb{E} \left[Y_{j,m}^2 \mathds{1}_{\left\{\left|Y_{j,m}\right|> \varepsilon s_n\right\}}\right] \longrightarrow 0, \ \ \text{as}\ n\to \infty.
		\end{aligned}
\end{equation}}
\textcolor{blue}{For $j=1,\ldots, n$, we have \begin{align*}
		\left|T_{j,m}\right|&\leq \sum \limits_{k=0}^{m-1} \left|\mathds{1}_{\left\{\frac{k}{m} < \hat{\phi}_{j,T}\leq \frac{k+1}{m}\right\}}\right|p_k(m-1,x)+  \sum \limits_{k=0}^{m-1} \left|G\left(\dfrac{k+1}{m}\right)- G\left(\dfrac{k}{m}\right)\right| p_k\left(m-1,x\right)\\ 
		& \leq \underset{0\leq k\leq m-1}{\max} p_k\left(m-1,x\right)+\underset{0\leq k\leq m-1}{\max}\left( G\left(\dfrac{k+1}{m}\right)- G\left(\dfrac{k}{m}\right) \right)\\
		&=O(m^{-1})+\left(\sum \limits_{k=0}^{m-1} p^2_k\left(m-1,x\right)\right)^{\frac{1}{2}}\\
		&=O(m^{-\frac{1}{4}}).
\end{align*}}
\textcolor{blue}{Hence, $\left|Y_{j,m}\right|=n^{-\frac{1}{2}}\left|T_{j,m}\right|=O(n^{-\frac{1}{2}}m^{-\frac{1}{4}}).$\\
	Further, for $x\in\left(0,1\right)$, we have  \begin{align*}
		s_n^2=n^{-1}\sum \limits_{j=1}^{n} \mathbb{E}\left(T^2_{j,m}\right)=\mathbb{E}\left(T^2_{1,m}\right)= Var\left(Z_{1,m}^T\right)
		=m^{-\frac{3}{2}}g(x)\psi(x)(1+o(1))
	\end{align*}
	Then \begin{align*}
		|Y_{j,m}|s_{n}^{-1}=O\left(n^{-\frac{1}{2}}m^{\frac{3}{4}}m^{-\frac{1}{4}}\right)=O\left(n^{-\frac{1}{2}}m^{\frac{1}{2}}\right){\longrightarrow }\ 0,
	\end{align*}
	whenever $mn^{-1}\rightarrow 0$ as $m,n \rightarrow \infty$. Under this condition, (\ref{lin}) holds and by Lindeberg-Feller’s central limit theorem, we obtain\begin{align*}
		s_{n}^{-1} \sum \limits_{j=1}^{n} Y_{j,m} \overset{\mathcal{L}} {\longrightarrow } \mathcal{N}\left(0,1\right)\ \ \text{as}\ n\to \infty.
\end{align*}}
\textcolor{blue}{It follows that $$n^{\frac{1}{2}}m^{-\frac{1}{4}} \left(\hat{f}_{n,m}(x)-\mathbb{E}\left(\hat{f}_{n,m}(x)\right)\right) \overset{\mathcal{L}} {\longrightarrow } \mathcal{N}\left(0,g(x)\psi(x)\right),\ \text{for}\ x\in\left(0,1\right).$$
	For the second assertion of Theorem \ref{Teo1}, we have for $x\in\left\{0,1\right\}$\\
	$$\left|Y_{j,m}\right|=n^{-\frac{1}{2}}\left|T_{j,m}\right|\leq 2n^{-\frac{1}{2}}\ \ 
	\text{and}\ \ s_n^2=m^{-1}g(x)+ O(m^{-1}).$$
	It follows that $|Y_{j,m}|s_{n}^{-1}=O\left(n^{-\frac{1}{2}}m^{\frac{1}{2}}\right)$ which vanishes when $m,n \rightarrow \infty$ such that $mn^{-1}\rightarrow 0$.
	Thus, by Lindeberg-Feller's central limit theorem, we get $$n^{\frac{1}{2}}m^{-\frac{1}{2}}\left[\hat{f}_{m,n}(x)- \mathbb{E}\left(\hat{f}_{m,n}(x)\right)\right] \overset{\mathcal{L}}{\longrightarrow} \mathcal{N}\left(0,g(x)\right),$$
	which concludes the proof of Theorem\ref{Teo1}.}
\textcolor{blue}{\subsection{Proof of Corollary \ref{co2}}
	In order to prove this theorem, we split the term $\hat{f}_{m,n}(x)-f(x)$ into a variance part $\hat{f}_{m,n}(x)- \mathbb{E}\left(\hat{f}_{m,n}(x)\right)$ and a bias part $\mathbb{E}\left(\hat{f}_{m,n}(x)\right)-f(x)$.
	In order to prove the first result of Theorem \ref{co2}, we write for all $x\in\left(0,1\right)$
	\begin{equation}
		\label{dec}
		n^{\frac{1}{2}} m^{-\frac{1}{4}} \left( \hat{f}_{m,n}(x)-f(x)\right)=n^{\frac{1}{2}} m^{-\frac{1}{4}} \left( \hat{f}_{m,n}(x)- \mathbb{E}\left(\hat{f}_{m,n}(x)\right)\right) + n^{\frac{1}{2}} m^{-\frac{1}{4}} \left(\mathbb{E}\left(\hat{f}_{m,n}(x)\right)-f(x)\right)
\end{equation}}
\textcolor{blue}{Under assumption $(A4)$, the bias term can be written for all $x\in[0,1]$ as follows
	\begin{equation}
		\label{biasg}
		\mathbb{E}\left(\hat{f}_{m,n}(x)\right)-f\left(x\right)=m\sum\limits_{k=0}^{m-1}\left(G\left(\frac{k+1}{m}\right)-G\left(\frac{k}{m}\right)\right)p_k\left(m-1,x\right)-f(x).\end{equation}}
\textcolor{blue}{
	By using the Taylor Young expansion, we obtain
	\[
	G\left(\frac{k+1}{m}\right)-G\left(\frac{k}{m}\right)=m^{-1}g\left(\frac{k}{m}\right)+O\left(m^{-1}\right),\]
	and \[g\left(\frac{k}{m}\right)=g(x)+ O\left(\left|\frac{k}{m}-x\right|\right),\]
	Substituting these two expressions in (\ref{biasg}), we get
	\[\mathbb{E}\left(\hat{f}_{m,n}(x)\right)-f\left(x\right)=g\left(x\right)-f\left(x\right)+O\left(m^{-\frac{1}{2}}\right)+O\left(1\right)\]}
\textcolor{blue}{Now, by replacing this term into (\ref{dec}), we get 
	\begin{align*}
		n^{\frac{1}{2}} m^{-\frac{1}{4}} \left( \hat{f}_{m,n}(x)-f(x)\right)=n^{\frac{1}{2}} m^{-\frac{1}{4}} \left( \hat{f}_{m,n}(x)- \mathbb{E}\left(\hat{f}_{m,n}(x)\right)\right) + n^{\frac{1}{2}} m^{-\frac{1}{4}}\left(g\left(x\right)-f\left(x\right)\right)+O\left(n^{\frac{1}{2}} m^{-\frac{1}{4}}\right).
	\end{align*}
	The first assertion of Theorem \ref{Teo1} ensures that as $m,n\ \to \infty$ such that $mn^{1}\to 0$, the variance term verifies \[n^{\frac{1}{2}} m^{-\frac{1}{4}} \left( \hat{f}_{m,n}(x)- \mathbb{E}\left(\hat{f}_{m,n}(x)\right)\right) \overset{\mathcal{L}}{\longrightarrow} \mathcal{N}\left(0,g(x)\psi(x)\right). \] 
	Moreover, if there exists a constant $c\geq0$ such that $n^{\frac{1}{2}} m^{-\frac{1}{4}}\to c$ as $m,n\to\infty$, the second term of the right-hand side of (\ref{dec}) converges in probability to $\delta=c\left(g\left(x\right)-f\left(x\right)\right)$.
	Then by Slutsky's lemma we get \begin{align*}
		n^{\frac{1}{2}} m^{-\frac{1}{4}} \left( \hat{f}_{m,n}(x)-f(x)\right) \overset{\mathcal{L}}{\longrightarrow} \mathcal{N}\left(\delta, g(x)\psi(x)\right).
	\end{align*}
	Using the same steps, we can prove the result stated in the second assertion of our theorem.}
	\textcolor{blue}{\subsection{Proof of Theorem \ref{t3}} 
		Let \(\tilde{f}_{m,n}\) be the Bernstein density estimator based on the unobserved random effects \(\phi_j\) as defined in (\ref{estb}) and denote $f_m(x):=\mathbb{E}\left(\tilde{f}_{m,n}(x)\right)$. Under these notations, we obtain
		\begin{equation}
			\label{norm}
			\mathbb{E} \left\|\hat{f}_{m,n}-f \right\|\leq \mathbb{E} \left\|\hat{f}_{m,n}-\tilde{f}_{m,n} \right\| + \mathbb{E} \left\|\tilde{f}_{m,n}- f_m \right\|+\left\|f_m-f \right\|.
		\end{equation}
		The first term in the right hand side of (\ref{norm}) is due to the estimation of $\phi_j$ by $\hat{\phi}_{j,T}$, the second one is a variance term and the last one corresponds to the bias term.
		In the sequel we aim to find an upper bound for each of these terms. To that end, we study each bound in a separate paragraph.}
	\textcolor{blue}{\paragraph{Upper bound for \(\mathbb{E} \left\|\hat{f}_{m,n}-\tilde{f}_{m,n} \right\|\):}
		\begin{equation}
			\begin{aligned}
				\label{term1n}
				\mathbb{E} \left\|\hat{f}_{m,n}-\tilde{f}_{m,n} \right\|&= \mathbb{E}\left[\sup_{x \in [0,1]} \left| \dfrac{m}{n} \sum_{k=0}^{m-1}\sum\limits_{j=1}^n \left(\mathds{1}_{\left\{\hat{\phi}_{j,T} \in\beta_k\right\}}-\mathds{1}_{\left\{\phi_j\in \mathcal{B}_k\right\}}\right) p_k(m-1, x) \right|\right]\\
				&=\mathbb{E} \left[m\sup_{x \in [0,1]} \left| \sum_{k=0}^{m-1} Y_k p_k(m-1, x) \right|\right],
			\end{aligned}
		\end{equation}
		where
		\[
		Y_k = \frac{1}{n} \sum_{j=1}^n \left( \mathds{1}_{\left\{\hat{\phi}_{j,T}\in \mathcal{B}_k\right\}} - \mathds{1}_{\left\{\phi_j\in \mathcal{B}_k\right\}} \right),\ \text{for all}\ 0\leq k \leq m-1.
		\]
		Taking into account the fact that Bernstein polynomials satisfy for all $x\in[0,1]$
		\begin{equation}
			\label{propbern}
			\sum_{k=0}^{m-1} p_k(m-1, x) = 1\ \text{and}\ \quad 0 \leq p_k(m-1, x) \leq 1, \ \text{for}\ 0\leq k\leq m-1,
		\end{equation}
		we can further bound (\ref{term1n}) as follows
		\begin{equation}
			\begin{aligned}
				\label{maj1}\mathbb{E} \left[m \sup_{x \in [0,1]} \left| \sum_{k=0}^{m-1} Y_k p_k(m-1, x) \right|\right]&\leq  m\mathbb{E} \left[ \sup_{x \in [0,1]} \sum_{k=0}^{m-1} \left|Y_k \right|p_k(m-1, x) \right]\nonumber\\
				&\leq m \mathbb{E} \left[ \max_{0\leq k \leq m-1}  \left|Y_k \right| \right],
			\end{aligned}
	\end{equation}}
	\textcolor{blue}{ which implies that
		\begin{equation}
			\label{finalmaj1}
			\mathbb{E} \left\|\hat{f}_{m,n}-\tilde{f}_{m,n}\right\|\leq m \mathbb{E} \left[ \max_{0\leq k \leq m-1}  \left|Y_k \right| \right].
		\end{equation}
		Our aim now is to control the right-hand side expectation term of \ref{finalmaj1}. 
		Since the difference \( \mathds{1}_{\hat{\phi}_{j,T} \in \beta_k} - \mathds{1}_{\phi_j \in \beta_k} \) is non-zero only if \( \hat{\phi}_{j,T} \) and \( \phi_j \) fall into different intervals, its expectation can be bounded, for $T>0$, $1\leq j\leq n$ and $0\leq k \leq m-1$, as follows
		\begin{eqnarray*}
			\mathbb{E} \left| \mathds{1}_{\left\{\hat{\phi}_{j,T} \in\beta_k\right\}} - \mathds{1}_{\left\{\phi_j\in \beta_k\right\}} \right| = \mathbb{P} \left( \hat{\phi}_{j,T} \notin \beta_k, \phi_j \in \beta_k \right) + \mathbb{P} \left( \hat{\phi}_{j,T} \in \beta_k, \phi_j \notin \beta_k \right)
			\leq \mathbb{P}\left( |\hat{\phi}_{j,T} - \phi_j| \geq \frac{1}{m} \right).
		\end{eqnarray*}
		Since $\hat{\phi}_{j,T}$ is an unbiased estimator for $\phi_j$, the difference \(\hat{\phi}_{j,T}-\phi_j\) is a centered random variable with variance giving by 
		\begin{equation}
			\label{varerreur}
			Var\left(\hat{\phi}_{j,T}-\phi_j\right)=\mathbb{E}\left[ (\hat{\phi}_{j,T} - \phi_j)^2 \right],\end{equation}
		and from equation (\ref{er}), we have for $T>0$ and $1\leq j\leq n$
		\[
		\mathbb{E}\left[ (\hat{\phi}_{j,T} - \phi_j)^2 \right]= \dfrac{\mathbb{E}\left[\left(\int_0^T J_2(t) dM_t^{H,j}\right)^2\right]}{\left(\int_0^T J_2^2(t) dw_t^H\right)^2},
		\]}
	\textcolor{blue}{where $w_t^H$, $J_2(t)$ and $M^{H,j}$ are as defined respectively in (\ref{2}), (\ref{4}) and (\ref{3}).
		Since the Molchan martingale $M^{H,j}$ is a continuous fractional martingale, with quadratic variation $\langle M^{H,j}_T \rangle=w_T^H$, it follows from the isometry of stochastic integrals that 
		\[\mathbb{E}\left[\left(\int_0^T J_2(t) dM_t^{H,j}\right)^2\right]=\int_0^T J_2^2(t)dw_t^H\]
		which implies that 
		\[\mathbb{E}\left[ (\hat{\phi}_{j,T} - \phi_j)^2 \right]=\left(\int_0^T J_2^2(t) dw_t^{H}\right)^{-1}.\]
		Assume that for all $s\in[0,T]$, there exists a positive constant $C$ such that $\dfrac{b(s)}{\sigma(s)}\geq C$. Then, \[
		\int_0^T J_2^2(t) d\,w_t^{H}=\int_0^T \left(\dfrac{d}{d\,w_t^H}\int_0^t k(t,s)\dfrac{b(s)}{\sigma(s)} ds\right)^2 d\,w_t^{H} \geq C^2 \int_0^T \left(\dfrac{d}{d\,w_t^H}\int_0^tk(t,s) ds\right)^2 d\,w_t^{H},
		\]
		where $k(t,s)$ is as defined in (\ref{1}). As shown in Example \ref{examp} the left hand-side integral is equal to \[\int_0^T \left(\dfrac{d}{d\,w_t^H}\int_0^tk(t,s) ds\right)^2 d\,w_t^{H}=\lambda_H^{-1} T^{2-2H},\]
		where $\lambda_H$ is as denoted in (\ref{2}). We then deduce that \[
		\mathbb{E}\left[ (\hat{\phi}_{j,T} - \phi_j)^2 \right]=\left(\int_0^T J_2^2(t) dw_t^{H}\right)^{-1}\leq \dfrac{\lambda_H}{C^2T^{2-2H}}.
		\]}
	\textcolor{blue}{Taking into account (\ref{varerreur}), we obtain for $T>0$ and $1\leq j\leq n$ \[Var\left(\hat{\phi}_{j,T}-\phi_j\right)\leq \dfrac{\lambda_H}{C^2T^{2-2H}}< \infty.\]\\
		Now by applying Chebyshev's inequality to the random variable \(\hat{\phi}_{j,T}-\phi_j\), we get for $m>0$
		\[
		\mathbb{P}\left( |\hat{\phi}_{j,T} - \phi_j| \geq \frac{1}{m} \right) \leq m^{2}\ \mathbb{E}\left[ (\hat{\phi}_{j,T} - \phi_j)^2 \right]\leq \dfrac{\lambda_H m^2}{C^2 T^{2 - 2H}},\]
		which in turn implies that for each $0\leq k\leq m-1$
		\[
		\mathbb{E} \left| \mathds{1}_{\left\{\hat{\phi}_{j,T} \in\beta_k\right\}} - \mathds{1}_{\left\{\phi_j\in \beta_k\right\}} \right| \leq \dfrac{\lambda_H m^2}{C^2 T^{2 - 2H}}.
		\]
		Therefore, for each $0\leq k\leq m-1$, we obtain\begin{equation*}
			\mathbb{E}\left(\left|Y_k\right|\right) \leq \dfrac{1}{n} \sum\limits_{j=1}^n \mathbb{E} \left| \mathds{1}_{\left\{\hat{\phi}_{j,T} \in\beta_k\right\}} - \mathds{1}_{\left\{\phi_j\in \beta_k\right\}} \right| \leq \dfrac{\lambda_H m^2}{C^2 T^{2 - 2H}}.
		\end{equation*}
		and to handle all $m$ intervals, a union bound implies that
		\[
		\mathbb{E}\left( \max_{0 \leq k \leq m-1} |Y_k| \right) \leq \sum\limits_{k=0}^{m-1}   \mathbb{E}\left(\left|Y_k\right|\right) \leq  \frac{\lambda_H m^3}{C^2 T^{2 - 2H}}.
		\]
		Substituting into (\ref{finalmaj1}), we get the following upper bound for the first term of (\ref{norm}) 
		\begin{equation}
			\mathbb{E} \left\|\hat{f}_{m,n}-\tilde{f}_{m,n} \right\| \leq  \frac{\lambda_H m^4}{C^2 T^{2 - 2H}}.
			\label{term1bound}
	\end{equation}}
	\textcolor{blue}{\paragraph{Upper bound for the variance term:}\ \\
		For the second term of (\ref{norm}), we have 
		\begin{equation}
			\begin{aligned}
				\label{term2n}
				\mathbb{E} \left\|\tilde{f}_{m,n}- f_m \right\|&= \mathbb{E}\left[\sup_{x \in [0,1]} \left| \dfrac{m}{n} \sum_{k=0}^{m-1}\sum\limits_{j=1}^n \left(\mathds{1}_{\left\{\phi_j \in\beta_k\right\}}-\mathbb{P}\left(\phi_j \in\beta_k\right)\right) p_k(m-1, x) \right|\right]\\
				&=\mathbb{E} \left[m\sup_{x \in [0,1]} \left| \sum_{k=0}^{m-1} H_k p_k(m-1, x) \right|\right],
			\end{aligned}   
		\end{equation}
		where
		\[
		H_k = \frac{1}{n} \sum_{j=1}^{n} \left(\mathds{1}_{\{\phi_j \in \beta_k\}} -  \mathbb{P}\left(\phi_j \in\beta_k\right)\right),\ \text{for all}\ 0\leq k \leq m-1.
		\]
		Using the properties stated in (\ref{propbern}), we can further bound (\ref{term2n}) as follows
		\begin{equation}
			\begin{aligned}
				\label{emax}
				\mathbb{E} \left[m \sup_{x \in [0,1]} \left| \sum_{k=0}^{m-1} H_k p_k(m-1, x) \right|\right]&\leq  \mathbb{E} \left[m \sup_{x \in [0,1]} \sum_{k=0}^{m-1} \left|H_k \right|p_k(m-1, x) \right]\\
				&\leq  m \mathbb{E} \left[ \max_{0\leq k \leq m-1}  \left|H_k \right| \right].
			\end{aligned}
		\end{equation}
		Thus, our task reduces to bounding $\displaystyle\mathbb{E} \left[ \max_{0\leq k \leq m-1}  \left|H_k \right| \right].$ Since $\displaystyle\max_{0\leq k \leq m-1}  \left|H_k \right|$ is a non-negative random variable, its expectation can be expressed as}
	\textcolor{blue}{
		\begin{equation}
			\label{espmax}
			\mathbb{E}\left[\max_{0 \leq k \leq m-1} |H_k|\right] = \int_0^\infty \mathbb{P}\left( \max_{0 \leq k \leq m-1} |H_k| \geq u\right) du.
	\end{equation}}
	\textcolor{blue}{Combining (\ref{term2n}), (\ref{emax}) and (\ref{espmax}), we get 
		\begin{equation}
			\label{finalmaj}
			\mathbb{E} \left\|\tilde{f}_{m,n}- f_m \right\|\leq m \int_0^\infty \mathbb{P}\left( \max_{0 \leq k \leq m-1} |H_k| \geq u\right) du.
		\end{equation}
		For each $0\leq k\leq m-1$, \( H_k \) is a centered random variable and its variance is expressed as \begin{equation}
			\begin{aligned}
				Var\left(H_k\right)= \mathbb{E}\left(H^2_k\right) &=\dfrac{1}{n^2} \mathbb{E}\left[\sum_{j=1}^{n} \left(\mathds{1}_{\{\phi_j \in \beta_k\}} -  \mathbb{P}\left(\phi_j \in\beta_k\right)\right)\right]^2\\
				&=\frac{1}{n^2}\mathbb{E} \left[\sum_{i=1}^{n}\sum\limits_{j=1}^n\left(\mathds{1}_{\left\{\phi_i \in\beta_k\right\}}-\mathbb{P}\left(\phi_i \in\beta_k\right)\right)\left(\mathds{1}_{\left\{\phi_j \in\beta_k\right\}}-\mathbb{P}\left(\phi_j \in\beta_k\right)\right)\right]\\
				&= \frac{1}{n^2} \sum_{i=1}^{n}\mathbb{E} \left[\left(\mathds{1}_{\left\{\phi_i \in\beta_k\right\}}-\mathbb{P}\left(\phi_i \in\beta_k\right)\right)^2\right]\\ &+ \frac{1}{n^2} \underset{i\neq j}{\sum}\mathbb{E} \left[\left(\mathds{1}_{\left\{\phi_i \in\beta_k\right\}}-\mathbb{P}\left(\phi_i \in\beta_k\right)\right)\left(\mathds{1}_{\left\{\phi_j \in\beta_k\right\}}-\mathbb{P}\left(\phi_j \in\beta_k\right)\right)\right]. 
				\label{firexp}
			\end{aligned}
		\end{equation}
		The first expectation of (\ref{firexp}) can be written as follows
		\begin{equation*}
			\begin{aligned}
				\mathbb{E}\left[\mathds{1}_{\left\{\phi_i \in\beta_k\right\}}-\mathbb{P}\left(\phi_i \in\beta_k\right)\right]^2&= \mathbb{P}\left(\phi_i \in\beta_k\right)  - 2\mathbb{P}\left(\phi_i \in\beta_k\right)^2 + \mathbb{P}\left(\phi_i \in\beta_k\right)^2.\\
				&=\mathbb{P}\left(\phi_i \in\beta_k\right)\left(1-\mathbb{P}\left(\phi_i \in\beta_k\right)\right),\label{secexp}
			\end{aligned}
		\end{equation*}
		and using the fact that the random effects $\phi_j$ are i.i.d, we obtain 
		\begin{equation*}
			\mathbb{E} \left[\left(\mathds{1}_{\left\{\phi_i \in\beta_k\right\}} - \mathbb{P}(\phi_i \in\beta_k)\right) \left(\mathds{1}_{\left\{\phi_j \in\beta_k\right\}} - \mathbb{P}(\phi_j \in\beta_k)\right)\right]=0.\label{thirexp}
		\end{equation*}
		Combining relations (\ref{firexp}), (\ref{secexp}) and (\ref{thirexp}), we obtain
		$$
		Var\left(H_k\right)=
		\frac{1}{n^2} \sum_{i=1}^{n}\mathbb{P}\left(\phi_i \in\beta_k\right)\left(1-\mathbb{P}\left(\phi_i \in\beta_k\right)\right)=\frac{1}{n}\mathbb{P}\left(\phi_1\in\beta_k\right)\left(1-\mathbb{P}\left(\phi_1\in\beta_k\right)\right).
		$$}
	\textcolor{blue}{Using the fact that $x(1-x)\leq \frac{1}{4}$ for $x\in [0,1]$ and since $0\leq \mathbb{P}\left(\phi_j \in\beta_k\right)\leq 1$ for all $0\leq k\leq m-1$, we deduce that
		$$\mathbb{P}\left(\phi_1 \in\beta_k\right)\left(1-\mathbb{P}\left(\phi_1 \in\beta_k\right)\right)\leq \dfrac{1}{4}.$$
		Then for each $0\leq k\leq m-1$, we obtain \(Var\left(H_k\right)\leq \dfrac{1}{4n}<\infty.\)\\
		Therefore, by Chebyshev's inequality, for any $u>0$ and $0\leq k\leq m-1$, we get\[
		\mathbb{P}(|H_k| \geq u) \leq \dfrac{Var\left(H_k\right)}{u^2}\leq \dfrac{1}{4nu^2}.
		\]
		Now, to handle all $m$ intervals, we apply the union bound
		\[
		\mathbb{P} \left( \max_{0 \leq k \leq m-1} |H_k| \geq u \right) \leq \sum\limits_{k=0}^{m-1} \mathbb{P}(|H_k| \geq u) \leq \dfrac{m}{4nu^2}.
		\]
		Replacing this bound into (\ref{espmax}), the expectation of the maximum can be bounded as follows
		\[
		\mathbb{E}\left[\max_{0 \leq k \leq m-1} |H_k|\right] = \int_0^\infty \mathbb{P}\left( \max_k |H_k| \geq u\right) du\leq 
		\int_0^\infty \min\left(1, \dfrac{m}{4nu^2} \right) du.
		\]
		Furthermore, the integral in the right-hand side can be split into two parts based on the value of $u$ that verifies $\dfrac{m}{4nu^2}=1$.\\ Denote $u_0$ the solution of this equation, it follows that
		\( 
		u_0 = \dfrac{1}{2}m^{\frac{1}{2}}n^{-\frac{1}{2}}.
		\)
		The integral then becomes
		\begin{equation*}
			\begin{aligned}
				\int_0^\infty \min\left(1, \dfrac{m}{4nu^2} \right) du&=\int_0^{u_0} 1 \, du + \int_{u_0}^{\infty} \dfrac{m}{4nu^2} \, du\\
				&= u_0 + \frac{m}{4n}\dfrac{1}{u_0}\\
				&= m^{\frac{1}{2}} n^{-\frac{1}{2}}.
			\end{aligned}
		\end{equation*}
		Finally, by replacing this integral into (\ref{finalmaj}), we obtain the following upper bound for the variance term
		\begin{equation}
			\label{variancebound}
			\mathbb{E} \left\|\tilde{f}_{m,n}- f_m \right\|\leq m^{\frac{3}{2}}n^{-\frac{1}{2}}.
	\end{equation}}
	\textcolor{blue}{\paragraph{Upper bound for the bias term:}\ \\
		For the last term in the right hand side of (\ref{norm}), we have for all $x\in[0,1]$
		\begin{equation}
			\label{terme3}
			\left\|f_m-f \right\|=\sup_{x \in [0,1]} \left|m \sum_{k=0}^{m-1} \left(F\left(\frac{k+1}{m}\right)-F\left(\frac{k}{m}\right)\right) p_k(m-1, x) - f(x) \right|.
		\end{equation}
		Since \( f \in C^2[0,1] \), applying Taylor-Young's expansion to second order for the function $F$, we obtain
		\begin{equation*}
			F\left(\frac{k+1}{m}\right) - F\left(\frac{k}{m}\right) = \frac{1}{m}f\left(\frac{k}{m}\right)+ \frac{1}{2m^2} f'\left(\frac{k}{m}\right) + O\left(m^{-3}\right).
		\end{equation*}
		By taking the sum overall $0\leq k\leq m-1$, we get
		\begin{eqnarray*}
			m \sum_{k=0}^{m-1} \left(F\left(\frac{k+1}{m}\right)-F\left(\frac{k}{m}\right)\right) p_k(m-1, x) &=& \sum_{k=0}^{m-1} f\left(\frac{k}{m}\right) p_k(m-1, x) + \frac{1}{2m} \sum_{k=0}^{m-1} f'\left(\frac{k}{m}\right) p_k(m-1, x)\\ &+& O\left(m^{-2}\right).
		\end{eqnarray*}
		Applying again Taylor-Young's expansion to second order for the function $f$, we can write 
		\[
		f\left(\frac{k}{m}\right) = f(x) + \left(\frac{k}{m} - x\right) f'(x) + \frac{1}{2} \left(\frac{k}{m} - x\right)^2 f''(x) + O\left(\left(\frac{k}{m} - x\right)^3\right).
		\]
		Then, \begin{eqnarray*}\sum_{k=0}^{m-1} f\left(\frac{k}{m}\right) p_k(m-1, x) &=& f(x) \sum_{k=0}^{m-1} p_k(m-1, x) + f'(x) \sum_{k=0}^{m-1} \left(\frac{k}{m} - x\right) p_k(m-1, x)\\ &+& \frac{1}{2} f''(x) \sum_{k=0}^{m-1} \left(\frac{k}{m} - x\right)^2 p_k(m-1, x) + O\left(m^{-3/2}\right).
		\end{eqnarray*}
		A last application of Taylor-Young's expansion for the function $f'$ allows us to write that \[
		f'\left(\frac{k}{m}\right) = f'(x) + \left(\frac{k}{m} - x\right) f''(x) + O\left(\left(\frac{k}{m} - x\right)^2\right),
		\]
		which implies that \begin{eqnarray*}
			\sum_{k=0}^{m-1} f'\left(\frac{k}{m}\right) p_k(m-1, x) = f'(x) \sum_{k=0}^{m-1} p_k(m-1, x) + f''(x) \sum_{k=0}^{m-1} \left(\frac{k}{m} - x\right) p_k(m-1, x) +O\left(m^{-1}\right).
		\end{eqnarray*}
	}
	\textcolor{blue}{Combining all the terms, we obtain
		\begin{eqnarray*}
			m \sum_{k=0}^{m-1} \left(F\left(\frac{k+1}{m}\right) - F\left(\frac{k}{m}\right)\right) p_k(m-1, x)&=&f(x)+ f'(x)\sum_{k=0}^{m-1}\left(\dfrac{k}{m}-x\right)p_k(m-1, x) +  \dfrac{1}{2m}f'(x)\\ &+&\dfrac{1}{2m}f''(x)\sum_{k=0}^{m-1}\left(\dfrac{k}{m}-x\right)^2 p_k(m-1, x) + O\left(m^{-1}\right).
		\end{eqnarray*}
		Taking into account that the Bernstein polynomials satisfy for all $x\in[0,1]$
		$$\sum_{k=0}^{m-1}\left(\dfrac{k}{m}-x\right)p_k(m-1, x)=0\ \  \text{and}\ \ \sum_{k=0}^{m-1}\left(\dfrac{k}{m}-x\right)^2 p_k(m-1, x)=\dfrac{x(1-x)}{m},$$
		it follows that $$ m \sum_{k=0}^{m-1} \left(F\left(\frac{k+1}{m}\right) - F\left(\frac{k}{m}\right)\right) p_k(m-1, x)=f(x) +\dfrac{x(1-x)}{2m}f''(x) + \dfrac{1}{2m}f'(x) + O\left(m^{-1}\right).$$
		Now taking the supremum over \([0,1]\) and using the fact that the function \(x\mapsto  x(1-x) \) attains its maximum on $[0,1]$ at \( x = \frac{1}{2} \), we obtain 
		\begin{equation}
			\begin{aligned}
				\label{biasbound}
				\sup_{x \in [0,1]} \left| m \sum_{k=0}^{m-1} \left(F\left(\frac{k+1}{m}\right) - F\left(\frac{k}{m}\right)\right) p_k(m-1, x)-f(x) \right|&=\sup_{x \in [0,1]} 
				\left| \frac{f'(x)}{2m}+\frac{f''(x)}{2m} x(1-x)  + O(m^{-1}) \right|\\&\leq \frac{\|f'\|}{2m} + \frac{\|f''\|}{8m}+ O(m^{-1}).
			\end{aligned}
		\end{equation}
		Since $f\in C^2[0,1]$,  $\|f'\|$ and $ \|f''\|$ are automatically finite due to the continuity of $f'$ and $f''$ on the compact interval $[0,1]$. Then the last term is well controlled and the bound holds.\\
		Combining the upper bound (\ref{term1bound}), (\ref{variancebound}) and (\ref{biasbound}), we obtain the desired result.
	}
	
	\appendix
	\setcounter{secnumdepth}{0}
	\section{Appendix}
	In what follows, we recall some limit theorems for continuous local martingales used in the study of the asymptotic behavior of the MLE of the random effects.
	First we recall a strong law of large numbers for continuous
	local martingales, see \cite{Lip01}.
	\begin{theor}[Liptser and Shiryaev (2001) ]\label{Liptser and Shiryaev}\ \\
		Let $(\Omega,\mathcal{F},(\mathcal{F})_{t\in \mathbb{R}_{+}},\mathbb{P})$ be a filtered probability space that meets the usual conditions. Let  $(M_{t})_{t \in \mathbb{R}_{+}}$ be a square-integrable continuous local martingale with respect to the filtration $(\mathcal{F})_{t\in \mathbb{R}_{+}}$ such that $\mathbb{P} \left( M_0 = 0\right)=1.$ Let $(\xi_t)_{t\in \mathbb{R}_{+}}$ be a progressively measurable process such that
		$$\mathbb{P}\left( \int_{0}^{t}  \xi_{u}^{2}\ d\left\langle M\right\rangle_u < \infty\right)=1, \quad t \in \mathbb{R}_{+},$$
		and
		\begin{align*}\label{e91}
			\int_{0}^{t}  \xi_{u}^{2}\, \ d \left\langle M\right\rangle_u \overset{a.s.}{\to} \infty,\quad as\;\; t\to \infty,
		\end{align*}
		where $\left( \left\langle M\right\rangle_t  \right)_{t \in \mathbb{R}_{+}} $ denotes the quadratic variation process of $M.$ Then
		\begin{align*}
			\dfrac{\int_{0}^{t}  \xi_{u}\, \ d M_u}{\int_{0}^{t}  \xi_{u}^{2}\, \ d\left\langle M\right\rangle_u}
			\overset{a.s.}{\to} 0,\quad as\;\; t\to \infty.
		\end{align*}
	\end{theor}
	In the next theorem, we present a generalized version for the central limit theorem that was investigated by Touati \cite{touati}. 
	To that aim, we consider a $d$-dimensional quasi-left continuous
	martingale $M=(M_{t})_{t \geq 0}$ locally square integrable, defined
	on a filtered space of probability
	$(\Omega,\mathcal{F},(\mathcal{F})_{t\geq 0},\mathbb{P})$ (see Jacod
	and Shiryaev \cite{jac03}) and we consider a deterministic $d\times
	d$ non-singular matrix process $V=(V_{t})_{t\geq 0}$. For $u \in
	\mathbb{R}^{d}$, we set
	\begin{align*}
		\Phi_{t}(u):=\exp \left(-\dfrac{1}{2}u^{*}\langle M^{c}
		\rangle_{t}\ ,u+\int_{0}^{t}\int_{\mathbb{R}^{d}}^{}(\exp\left\{i\langle
		u,x\rangle \right\}-1-i\langle u,x\rangle)\nu^{M}(d s, d x)\right),
	\end{align*}
	where $\left( \langle M^{c}\rangle_t\right) _{t \geq 0}$ is the quadratic variation of the continuous local martingale $M^{c}$ and $\nu^{M}$ denotes the compensator of the random measure
	associated to the jumps of $M$.\\
	In this Theorem and under an assumption  $(\mathcal{H})$ (see Theorem \ref{touati1} below) with respect to characteristic function $\phi$ given by \label{Fcarcteristique} instead of the {\it classic} Lindeberg condition, Touati established a  generalized central limits theorem for the  martingale $M$.\\
	Remember that  the classic CLT is obtained under assumption
	\begin{align*}
		(\mathcal{H}_{1})\quad V_{t}^{-1}\langle M \rangle_{t}(V_{t}^{*})^{-1}\overset{a.s.}{\longrightarrow}  C, \quad(t\longrightarrow\infty),
	\end{align*}
	and the Lindeberg condition
	\begin{align*}
		(\mathcal{H'})\quad \forall \delta > 0,
		\int_{\mathbb{R}^{d}}\int_{0}^{1}\parallel V_{t}^{-1}x\parallel^{2}\mathds{1}_{\left\{\parallel V_{t}^{-1}x\parallel t > \delta\right\}}\nu^{M}(ds,dx)\overset{a.s.}{\longrightarrow} 0,\quad(t\longrightarrow\infty),
	\end{align*}
	which imply assumption $(\mathcal{H})$ with
	$$\eta=C^{\frac{1}{2}} \quad \mbox{and}\quad \Phi_{\infty}(\eta,u)=\exp \left(-\dfrac{1}{2}u^{*}C u \right).$$
	Let us now state the theorem.
	\begin{theor}\label{touati1}{\empty}(GCLT)\ \\ Let $M=(M_{t})_{t\in \mathbb{R}_{+}}$ be a d-dimensional quasi-left
		continuous local martingale with $M_{0}=0$ and $V=(V_{t})_{t\in \mathbb{R}_{+}}$ a
		deterministic family of non-singular matrices. We define a
		probability $\mathcal{Q}$ on the space
		$\mathrm{C}(\mathcal{X},\mathbb{R}^{d})$ of continuous functions
		from $\mathcal{X}$ to $\mathbb{R}^{d}$ (where $\mathcal{X}$
		indicates a vector space of finite dimension). If
		the couple $(M,V)$ satisfies the following assumption\\\\
		$\mathcal{(H)}\quad \left\{
		\begin{array}{lll}
			\Phi_{t}((V_{t}^{T})^{-1}u) \overset{a.s.}{\longrightarrow}  \Phi_{\infty}(\eta,u),\quad \text{as} \quad t\to \infty,& &\\\\
			\Phi_{\infty}(\eta,u)\neq 0 \quad a.s.,&&
		\end{array}
		\right.$\vspace{0.5cm} \\where $\eta$ denotes a r.v., possibly
		degenerated taking values in $\mathcal{X}$ and
		$$\Phi_{\infty}(z,u)=\int_{\mathbb{R}^{d}}^{}\exp\{i< u,\xi>\}\pi(z,d\xi),\quad (z,u)\in \mathcal{X}\times \mathbb{R}^{d},$$ denotes the Fourier transform of the
		one-dimensional conditionals laws $(\pi(x,.),x\in\mathcal{X})$ of
		the probability $\mathcal{Q}$. Then
		$$(GCLT)~~~~~Z_{t}:=V_{t}^{-1}M_{t}\overset{\mathcal{L}}{\longrightarrow}
		Z_{\infty}:=\Sigma(\eta),\quad(t\to \infty),$$ in a
		stable manner where $(\Sigma(z),z\in\mathcal{X})$ is a $\mathcal{Q}$ law process independent of the r.v. $\eta$.
	\end{theor}It is important to note that proofs are easier to handle under $(\mathcal{H'})$ than under $(\mathcal{H})$.
	
	\newpage
	\bibliographystyle{plain}

\end{document}